\documentclass[11pt,reqno]{amsart}
\usepackage{graphicx}

\usepackage{amssymb}
\usepackage{cite}
\usepackage{amsmath}
\usepackage{latexsym}
\usepackage{amscd}
\usepackage{amsthm}
\usepackage{mathrsfs}
\usepackage{url}

\usepackage[backref=page,linktocpage=true,colorlinks,citecolor=magenta,linkcolor=blue,urlcolor=magenta]{hyperref}
\numberwithin{equation}{section}
\newtheorem{thm}{Theorem}
\newtheorem{corr}[thm]{Corollary}
\newtheorem{lem}[thm]{Lemma}
\newtheorem{prop}[thm]{Proposition}

\theoremstyle{definition}
\newtheorem{defn}{Definition}[section]
\theoremstyle{remark}
\newtheorem{rem}[thm]{Remark}%[section]
\newtheorem*{ack}{Acknowledgment}
\def\R{\mathbb R}

\def\SS{\mathbb S}
\def\W{\mathcal{W}}
\def\PW{\partial\mathcal{W}}

\def\f{\frac}

\def\pt{\partial}

\begin{document}
\title[Sharp bounds for the anisotropic $p$-capacity]{Sharp bounds for the anisotropic $p$-capacity of Euclidean compact sets}
%\author[H.~Li]{Haizhong~Li}
%\address{Department of Mathematical Sciences, Tsinghua University, Beijing 100084, P.~R.~China}
%\email{\href{mailto:lihz@tsinghua.edu.cn}{lihz@tsinghua.edu.cn}}
\author[R.~Li]{Ruixuan~Li}
\address{Department of Mathematical Sciences, Tsinghua University, Beijing 100084, P.~R.~China}
\email{\href{mailto:li-rx18@mails.tsinghua.edu.cn}{li-rx18@mails.tsinghua.edu.cn}}
\author[C.~Xiong]{Changwei~Xiong}
\address{College of Mathematics, Sichuan University, Chengdu 610065,  P.~R.~China}
\email{\href{mailto:changwei.xiong@scu.edu.cn}{changwei.xiong@scu.edu.cn}}
\date{\today}
%\thanks{This research was supported .}

\maketitle

\begin{abstract}
We prove various sharp bounds for the anisotropic $p$-capacity $\mathrm{Cap}_{F,p}(K)$ ($1<p<n$) of compact sets $K$ in the Euclidean space $\mathbb{R}^n$ ($n\geq 3$). For example, using the inverse anisotropic mean curvature flow (IAMCF), we get an upper bound of Szeg\"{o} type (1931) for $\mathrm{Cap}_{F,p}(K)$ when $\partial K$ is a smooth, star-shaped and $F$-mean convex hypersurface in $\mathbb{R}^n$ ($n\geq 3$). Moreover, for such a surface $\partial K$ in $\mathbb{R}^3$, by introducing the anisotropic Hawking mass and studying its monotonicity property along IAMCF, we obtain an upper bound of Bray--Miao type (2008) for $\mathrm{Cap}_{F,p}(K)$.
\end{abstract}

\section{Introduction}\label{sec1}

%Given $\Omega\subset \R^n$ ($n\geq 3$) a bounded open set, we define its Finsler $F$-capacity as
%\begin{equation}
%\mathcal{C}_F:=Cap_F(\Omega)=\inf\{\int_{\R^n}F(Dv)^2dx:v\in C_0^\infty (\R^n), \; v\geq 1 \text{ in } \Omega\}.
%\end{equation}
%Equivalently, we can define the anisotropic electrostatic capacity of $\Omega$ as $\mathcal{C}_F:=\int_U (F(\nabla u))^2 dx$ with $u\in E^1(U)$ being the unique solution of the PDE ($U:=\R^n\setminus \overline{\Omega}$)
%\begin{equation*}
%\begin{cases}
%\Delta_F u=0,\text{ in }U,\\
%u=1, \text{ on }\pt U,\\
%u\rightarrow 0,\text{ as } x\to \infty.
%\end{cases}
%\end{equation*}
%Some basic facts on the anisotropic capacity can be found in \cite{BCS16}.

The capacity problem is one of the most extensively-investigated topics in the potential theory, the mathematical physics, the partial differential equations, the convex geometry and other fields. The classical (electrostatic) capacity of a compact set $K$ in the Euclidean space $\R^3$ admits the physical interpretation that, it represents the maximal charge that can be put on $K$ while the electrical potential of the vector field created by this charge is no greater than one. See e.g. \cite{PS51,HKM06,Maz11,Shu,LL01} for some background information. In this paper we are concerned with sharp bounds of the anisotropic $p$-capacity for compact sets in the Euclidean space in terms of their various geometric quantities.

Let $F\in C^\infty(\R^n\setminus \{0\})\cap C(\R^n)$ be a Minkowski norm on $\R^n$ and $K \subset \mathbb{R}^n$ be a compact set. Throughout the paper we consider $n\geq 3$. For $1<p<n$, the \emph{anisotropic $p$-capacity} of $K$ is defined as
\begin{align}
\mathrm{Cap}_{F,p}(K) = \inf \{ \int_{\R^n} F^p(Dv) dx: v \in C_c^{\infty}(\R^n), v \geq 1\text{ on }K \},
\end{align}
where $C_c^\infty(\R^n)$ is the set of smooth functions with compact support in $\R^n$. See Section~\ref{sec2.1} for more details on the anisotropic $p$-capacity. In particular, by (2) in Proposition~\ref{prop-property}, $\mathrm{Cap}_{F,p}(K)=\mathrm{Cap}_{F,p}(\pt K)$ for a compact set. In view of this property some works in the literature are only concerned with the capacity for hypersurfaces that are the boundaries of compact sets.

Since the bounds we obtain for $\mathrm{Cap}_{F,p}(K)$ involve some anisotropic geometric quantities, let us first introduce them. Let $M$ be an immersed oriented hypersurface in $\R^n$ with the unit normal vector field $\nu$ and $d\mu$ be the area element of the induced metric on $M$. We can define the anisotropic area element on $M$ as $d\mu_F:=F(\nu)d\mu$ and the anisotropic area of $M$ as
\begin{align*}
|M|_F:=\int_M d\mu_F=\int_M F(\nu)d\mu.
\end{align*}
Meanwhile we define the anisotropic unit normal $\nu_F$ along $M$ as
\begin{align*}
\nu_F:=DF(\nu).
\end{align*}
We can check that $\nu_F$ lies on the boundary $\pt \mathcal{W}$ (named the Wulff shape) of the so-called Wulff ball $\mathcal{W}$, which is by definition the set
\begin{align*}
\mathcal{W}:=\{x\in \R^n:F^0(x)<1\}.
\end{align*}
Here $F^0$ is the dual norm of $F$, given as
\begin{align*}
F^0(x):=\sup_{\xi\neq 0}\frac{\langle x,\xi\rangle}{F(\xi)}.
\end{align*}
Moreover, for any point $x\in M$, the tangent hyperplane $T_xM$ is parallel to the tangent hyperplane $T_{\nu_F(x)}\pt \mathcal{W}$. So we may define the anisotropic Weingarten map of $M$ as
\begin{align*}
d\nu_F:T_xM\rightarrow T_{\nu_F(x)}\pt \mathcal{W}.
\end{align*}
The anisotropic Weingarten map $d\nu_F$ has $n-1$ real eigenvalues $\kappa^F_i$, $1\leq i\leq n-1$, which are called the anisotropic principal curvatures of the hypersurface. From them, for $1\leq k\leq n-1$ we can define the $k$th anisotropic mean curvature $\sigma_k(\kappa^F)$ of the hypersurface as
\begin{align*}
\sigma_k(\kappa^F):=\sum_{1\leq i_1<\cdots<i_k\leq n-1}\kappa^F_{i_1}\cdots \kappa^F_{i_k}.
\end{align*}
We make a convention $\sigma_0(\kappa^F)=1$. The special cases $H_F:=\sigma_1(\kappa^F)$ and $K_F:=\sigma_{n-1}(\kappa^F)$ are called the anisotropic mean curvature and the anisotropic Gaussian curvature of $M$, respectively. The hypersurface $M$ is called \emph{$F$-mean convex} if $H_F>0$ on the hypersurface.
\begin{rem}
For the isotropic case $F(\xi)=|\xi|$, all the geometric quantities above reduce to the ordinary ones in the Euclidean space.
\end{rem}
\begin{rem}
The Euclidean space $\R^n$ equipped with a Minkowski norm $F$ as the anisotropy makes a nice model for applications and produces fruitful results in the geometry and analysis. For instance, it may model an anisotropic medium where the growth of crystals, the noise-removal procedures in the digital image processing, the crystalline fracture theory, etc. can be studied. See e.g. \cite{BC18,BCS16} and references therein for some introductions.
\end{rem}

Now we are ready to state our main results in this paper. First we aim at the following sharp upper bounds for $\mathrm{Cap}_{F, p}(K)$.

\begin{thm}\label{thm1}
Let $K \subset \mathbb{R}^n$ ($n\geq 3$) be a compact set with non-empty interior. Suppose its boundary $\pt K$ is smooth, star-shaped and $F$-mean convex.

If $2\le p < n$, then
\begin{align}\label{eq-1}
\mathrm{Cap}_{F, p}(K) \leq \left(\frac{(p-1)(n-1)}{n-p}\right)^{1-p} \int_{\partial K} H_F^{p-1}\,d\mu_F.
\end{align}

If $1<p\leq 2\leq q<n$, then
\begin{align}\label{eq-2}
\mathrm{Cap}_{F, p}(K) \le \left( \frac{(p-1)(n-1)}{n-p} \right)^{1-p} \left(\int_{\partial K} H_F^{q-1}\,d\mu_F\right)^{\frac{p-1}{q-1}} |\partial K|_F^{\frac{q-p}{q-1}} .
\end{align}
Moreover, the equality holds in (\ref{eq-1}) or (\ref{eq-2}) if and only if $\pt K$ is a translated scaled Wulff shape.
\end{thm}
\begin{rem}
For the isotropic case $F(\xi)=|\xi|$, the first result in the spirit of Theorem~\ref{thm1} may go back to Szeg\"{o}'s \cite{Sze31} (see also \cite{PS51}), where he proved the upper bound in the case $n=3$ and $p=2$ for smooth convex compact sets. In \cite{FS14}, Freire and Schwartz obtained the bound in Theorem~\ref{thm1} for the isotropic case with $n\geq 3$ and $p=2$ for smooth compact sets with mean-convex and outer-minimizing boundary. Last, Xiao in \cite{Xiao17} proved exactly Theorem~\ref{thm1} for the isotropic case.
\end{rem}
\begin{rem}
Recently, Xia and Yin in \cite{XY20} considered the anisotropic case and derived an interesting related upper bound of $\mathrm{Cap}_{F, p}(K)$ for smooth compact connected sets in $\R^n$ ($n\geq 3$). Their result \cite{XY20} and our Theorem~\ref{thm1} are not mutually inclusive. Besides, we use a different method from theirs.
\end{rem}

Before stating our next result, we define the anisotropic Hawking mass $m_H^F(\Sigma)$ for an immersed smooth oriented compact surface $\Sigma\subset \R^3$ without boundary as
\begin{align}\label{eq-Hawking}
	m_H^F(\Sigma) := \sqrt{\frac{|\Sigma|_F}{4 |\PW|_F}} \left(1 - \frac{\int_{\Sigma} H_F^2\,d\mu_F}{4|\PW|_F}\right).
\end{align}
Note that if $\Sigma$ is embedded, then $m_H^F(\Sigma)\leq 0$ with the equality if and only if $\Sigma$ is a translated scaled Wulff shape; see Corollary~\ref{corr-non-positive} below.

The classical (isotropic) Hawking mass was introduced by S.~W.~Hawking in \cite{Haw68} and is one of the most important concepts in the mathematical physics. In \cite{Gero73} Geroch discovered the significant monotonicity property of the Hawking mass along the smooth inverse mean curvature flow (IMCF) in the manifold setting. This monotonicity property was extended by Huisken and Ilmanen \cite{HI01} to the weak IMCF to prove the famous Riemannian Penrose inequality. Later Bray and Miao \cite{BM08} applied the monotonicity property of the Hawking mass from \cite{HI01} to obtain a new sharp upper bound for the capacity in the manifold setting. See \cite{Xiao16} for the $p$-capacity generalization of \cite{BM08}.

Here we find an application of the anisotropic analogue of the Hawking mass in the Euclidean space as in the second main result of our paper, and we hope this new concept would be of use in other problems in the future. Our result is the anisotropic version of some results from \cite{BM08} and \cite{Xiao16} in $\R^3$.
\begin{thm}\label{thm2}
Let $K \subset \mathbb{R}^3$ be a compact set with non-empty interior. Suppose its boundary $\Sigma=\pt K$ is smooth, star-shaped and $F$-mean convex. Let $1<p<3$.

If the anisotropic Hawking mass $m_H^F (\Sigma) = 0$, then $K=r_0\overline{\mathcal{W}}+x_0$ for some $r_0>0$ and $x_0\in \R^3$, and
\begin{align}\label{case-zero}
\mathrm{Cap}_{F, p}\left(K\right) = \left(\frac{3-p}{p-1}\right)^{p-1} |\PW|_F \; r_0^{3-p}.
\end{align}

If $m_H^F (\Sigma) < 0$, then
\begin{align}
\mathrm{Cap}_{F, p}(K) &<  \left( \frac{3-p}{p-1} \right)^{p-1} |\PW|_F^{\frac{p-1}{2}} |\Sigma|_F^{\frac{3-p}{2}}\nonumber\\
 &\times \left( \frac{\int_{\Sigma} H_F^2\,d\mu_F}{4|\PW|_F}  - 1\right)^{3-p}  \theta^{1-p},\label{case-negative}
\end{align}
where $\theta$ is defined as
\begin{align*}
&\theta:= \int_{0}^{\left( \frac{\int_{\Sigma} H_F^2\,d\mu_F}{4|\PW|_F}  -1 \right)^{\frac{3-p}{p-1}}} \left(1 +  r^{\frac{p-1}{3-p}}\right)^{-\frac{1}{2}}\,dr.
\end{align*}
%The equality in \eqref{case-zero} holds if and only if $\Sigma$ is a Wulff Shape.
\end{thm}
\begin{rem}
The surface $\Sigma$ in Theorem~\ref{thm2} is embedded. So by the remark after \eqref{eq-Hawking}, we know $m_H^F(\Sigma)\leq 0$. In addition, the bound \eqref{case-negative} is sharp in the sense that as $m_H^F(\Sigma)\rightarrow 0^-$, it reduces to the equality case \eqref{case-zero}.
\end{rem}

When $p=2$, we obtain the following partial anisotropic generalization of a result in \cite{BM08}.
\begin{corr}\label{capf}
Let $K \subset \mathbb{R}^3$ be a compact set with non-empty interior. Suppose its boundary $\Sigma=\pt K$ is smooth, star-shaped and $F$-mean convex. Then we have
\begin{align}
\mathrm{Cap}_{F,2}\left(K\right) \le \frac{1}{2} \sqrt{|\mathcal{\pt W}|_F |\Sigma|_F} \left(1 + \sqrt{\frac{1}{4 |\mathcal{\pt W}|_F} \int_{\Sigma} H_F^2\,d\mu_F}\right).
\end{align}
The equality holds if and only if $\Sigma$ is a translated scaled Wulff Shape.
\end{corr}

Third, we obtain the following result, which is motivated by Sections~3.4 and 3.5 of P\'{o}lya and Szeg\"{o}'s book \cite{PS51}.
\begin{thm}\label{thm3}
Let $K\subset \R^n$ ($n\geq 3$) be a compact set with non-empty interior and with smooth boundary and let $1<p<n$.
	\begin{enumerate}
		\item Assume that $\pt K$ is convex. Then we have
		\begin{align}\label{bound2}
			\mathrm{Cap}_{F, p}(K) \leq \left(\int_0^\infty \left(\sum_{i=0}^{n-1}\int_{\pt K} \sigma_i(\kappa^F) d\mu_F \cdot t^i\right)^{1/(1-p)}dt\right)^{1-p} ,
		\end{align}
		where $\sigma_i(\kappa^F)$ ($0\leq i\leq n-1$) is the $i$th anisotropic mean curvature of $\pt K$. Moreover, the equality holds if and only if $\pt K$ is a translated scaled Wulff shape.
		\item Assume that $\pt K$ is star-shaped with respect to the origin. Then we have
		\begin{equation}\label{bound3}
			\mathrm{Cap}_{F, p}(K)\leq   \left(\frac{n-p}{p-1}\right)^{p-1}\int_{\pt K} h_F^{1-p}d\mu_F,
		\end{equation}
		where $h_F=\langle X,\nu\rangle /F(\nu)$ is the anisotropic support function of the hypersurface $\pt K$. Moreover, the equality holds if and only if $\pt K$ is a scaled Wulff shape centered at the origin.
	\end{enumerate}
\end{thm}
\begin{rem}
		For Case (1), when $n=3$ and $p=2$, we obtain
		\begin{align*}
			&\int_0^\infty \left(\sum_{i=0}^{n-1}\int_{\pt K} \sigma_i(\kappa^F) d\mu_F \cdot t^i\right)^{1/(1-p)}dt\\
=&\int_0^\infty \frac{1}{|{\pt K}|_F+\int_{\pt K} \sigma_1(\kappa^F) d\mu_F\cdot t+|\PW|_F\cdot t^2}dt\\
			=&\frac{1}{\int_{\pt K}\sigma_1(\kappa^F)d\mu_F\cdot \varepsilon}\log\frac{1+\varepsilon}{1-\varepsilon},
		\end{align*}
		where
		\begin{equation*}
			\varepsilon=\sqrt{1-\frac{4|\PW|_F |{\pt K}|_F}{(\int_{\pt K} \sigma_1(\kappa^F)d\mu_F)^2}}\in [0,1).
		\end{equation*}
		Consequently,
		\begin{align*}
			\mathrm{Cap}_{F,2}(K)\leq  \int_{\pt K}\sigma_1(\kappa^F) d\mu_F\frac{\varepsilon}{\log((1+\varepsilon)/(1-\varepsilon))}\leq \frac{\int_{\pt K}\sigma_1(\kappa^F) d\mu_F}{2}.
		\end{align*}
This is the anisotropic analogue of the original result due to Szeg\"{o} \cite{Sze31} in 1931.
		%For general $n\geq 3$, we have
%		\begin{align*}
%			T(t)&=\sum_{i=0}^{n-1}\int_M \sigma_i d\mu \cdot t^i\leq \sum_{i=0}^{n-1}\omega_{n-1}C_{n-1}^i(\frac{\int_M \sigma_{n-2}d\mu}{(n-1)\omega_{n-1}})^{n-1-i}t^i\\
%			&=\omega_{n-1}(t+\frac{\int_M \sigma_{n-2}d\mu}{(n-1)\omega_{n-1}})^{n-1},
%		\end{align*}
%		where $\omega_{n-1}=|\SS^{n-1}|$. Then we can derive
%		\begin{equation*}
%			\xi_1\leq  (n-2)\omega_{n-1}\left(\frac{\int_M \sigma_{n-2}d\mu}{(n-1)\omega_{n-1}}\right)^{n-2},
%		\end{equation*}
%		which is weaker than \eqref{bound2} in Theorem~\ref{thm1}.
	\end{rem}
\begin{rem}
In \cite{PS51}, P\'{o}lya and Szeg\"{o} did not consider the equality case of \eqref{bound2} or \eqref{bound3}. Here for the equality case of \eqref{bound2}, we need the rigidity result Theorem~1.2 of \cite{BC18}; while for that of \eqref{bound3}, we are inspired by Theorem~1.2 of \cite{BCS18}.
\end{rem}

For the proofs of Theorems~\ref{thm1}, \ref{thm2} and \ref{thm3}, the idea may originate from Szeg\"{o}'s work \cite{Sze31}. Let $U=\R^n\setminus K$. In each case we use a smooth family $\{M_t\}_{t\geq 0}$ with $M_0=\pt K$ of hypersurfaces to foliate $\overline{U}$ and construct a suitable test function $f(x)$ with level sets being $M_t$. Then the problem to find an upper bound of $\mathrm{Cap}_{F, p}(K)$ is reduced to the estimate for certain geometric quantity on the hypersurface $M_t$. Employing different techniques for such an estimate gives rise to different results as in Theorems~\ref{thm1}, \ref{thm2} and \ref{thm3}. More precisely, for Theorem~\ref{thm1}, we need the inverse anisotropic mean curvature flow (IAMCF) \cite{Xia17}. For Theorem~\ref{thm2}, we rely further on the monotonicity property of the anisotropic Hawking mass along the IAMCF. And for Theorem~\ref{thm3}, we construct a natural flow in each case.

Last, we derive the following sharp lower bound of $\mathrm{Cap}_{F, p}(K)$, which is a generalization of \cite[Theorem~2.1]{Xiao17}.
\begin{thm}\label{thm4}
Let $K \subset \mathbb{R}^n$ ($n\geq 3$) be a compact convex set with non-empty interior and $1<p< n$. Then
\begin{align}\label{eq-lower-bound}
\frac{p(n-1)}{n(n-p)}|\partial K|_F^{n/(n-1)}& |\pt \W|_F^{1/(1-n)} \geq  |K|+ \frac{|\partial K|_F^{p/(p-1)}}{\mathrm{Cap}_{F, p}(K)^{1/(p-1)}}.
\end{align}
The equality holds if and only if $\pt K$ is a translated scaled Wulff shape.
\end{thm}
\begin{rem}
For the set $K$ in Theorem~\ref{thm4}, it is well-known that the outward unit normal $\nu$ along $\pt K$ is well-defined almost everywhere; see e.g. \cite{Sch14}. So $|\pt K|_F$ is well-defined. In addition, see \eqref{eq-mixed-volume} below for its expression as a mixed volume.
\end{rem}
\begin{rem}
If $K \subset \mathbb{R}^n$ is an arbitrary compact convex set (possibly with empty interior), we may consider the compact convex set $K_\varepsilon:=\{z\in \R^n:z=x+\varepsilon y,x\in K,y\in \overline{B_1}\}$ for $\varepsilon>0$. Applying first Theorem~\ref{thm4} to $K_\varepsilon$ and then taking $\varepsilon\to 0^+$, we get the result for $K$ itself (in light of Proposition~\ref{prop-property} below).
\end{rem}
\begin{rem}
For the isotropic case $F(\xi)=|\xi|$, Xiao \cite{Xiao17} proved the inequality~\eqref{eq-lower-bound}, which leads to a crucial step towards the P\'{o}lya--Szeg\"{o} conjecture. See e.g. \cite{Xiao17,PS51,Pol47} for information on this important conjecture. On the other hand, in \cite{Xiao17}, Xiao did not handle the ``only if'' part of the equality case. Here we are able to do it because we have a key observation \eqref{eq-key-observation} in Section~\ref{sec6}.
\end{rem}
The proof of Theorem~\ref{thm4} follows closely that of \cite[Theorem~2.1]{Xiao17}. In the proof we study the super-level set $K_t:=\{x|u(x)\geq t\}$ ($0<t<1$) of the anisotropic $p$-capacitary potential $u(x)$ associated with $K$. The main tools we use include the relationship between $\mathrm{Cap}_{F, p}(K)$ and $\mathrm{Cap}_{F, p}(K_t)$ and the anisotropic isocapacitary inequality. See Section~\ref{sec6} for details.

Finally, we like to mention that some related works on the estimates of the capacity can be found in \cite{AFM20,AM15,AM20,PS51,Xiao17,XY20,Xiao16,FS14,BM08,Maz11,HKM06} and references therein. Meanwhile, as a concluding remark for the Introduction, it is worth highlighting that although our methods are from the literature, our results are consequences of highly non-trivial works \cite{And01,Xia17,BC18,AGHLV17,XY20} etc. and we have got some completely new results, e.g., the equality cases in Theorems~\ref{thm3} and \ref{thm4}. Besides, to the best of our knowledge, so far there have been few estimates on the anisotropic $p$-capacity in addition to the classical ones in \cite{Maz11,HKM06} and the recent ones in \cite{XY20}. We hope our results would be a nice stimulation in the field of estimates on the anisotropic $p$-capacity.

The paper is structured as follows. In Section~\ref{sec2} we review some basic facts on the anisotropic $p$-capacity, the anisotropic geometry of hypersurfaces in the Euclidean space and the inverse anisotropic mean curvature flow. In Section~\ref{sec3} we first introduce a general approach to attack Theorems~\ref{thm1}, \ref{thm2}, \ref{thm3}, and then prove Theorem~\ref{thm1}. In the next Sections \ref{sec4} and \ref{sec5} we prove Theorems~\ref{thm2} and \ref{thm3} respectively along this general approach. In the final Section~\ref{sec6} we prove Theorem~\ref{thm4} following the method in \cite{Xiao17}. Throughout the paper, the Einstein convention for the summation of indices is used unless otherwise stated, and we usually use $M$ to denote a hypersurface in $\R^n$ while $\Sigma$ to denote a surface in $\R^3$.

\begin{ack}
The authors wish to thank Prof.~Haizhong~Li for enlightening discussions and helpful comments. The second-named author would also like to thank Prof.~Ben~Andrews for his interest in the work and Prof.~Chao~Xia for the communication on some anisotropic problems. The work was started when the second-named author was a postdoctoral fellow at the Mathematical Sciences Institute, Australian National University, to which he is grateful for providing a wonderful working environment. This research was supported by the funding (no.~1082204112549) from Sichuan University.
\end{ack}

\section{Preliminaries}\label{sec2}

This section is devoted to a brief overview of some preliminary materials required in this paper, including the anisotropic $p$-capacity for sets in the Euclidean space, the classical anisotropic geometry of Euclidean hypersurfaces in the differential geometry, and the relatively new anisotropic geometry of Euclidean hypersurfaces in the geometric analysis together with the resulting inverse anisotropic mean curvature flow.

\subsection{Anisotropic $p$-capacity}\label{sec2.1}
For this subsection nice references include \cite[Section~2.2]{Maz11} and \cite[Chapters~2 and 5]{HKM06}. First we introduce the Minkowski norm on $\R^n$.
\begin{defn}
A function $F \in C^{\infty}(\mathbb{R}^{n} \setminus \{0\})\cap C(\R^n)$ is called a \emph{Minkowski norm} if
\begin{enumerate}
  \item $F$ is a convex, even, $1$-homogeneous function, and $F(\xi) > 0$ if $\xi\neq 0$;
  \item $F$ satisfies the uniformly elliptic condition, i.e., $\mathrm{Hess}_{\R^{n}}(F^2)$ is positive definite in $\R^{n}\setminus \{0\}$.
\end{enumerate}
\end{defn}

Let $K\subset \R^n$ be a compact set. For $n\geq 3$ and $1<p<n$, the \emph{anisotropic $p$-capacity} of $K$ is defined as
\begin{align*}
\mathrm{Cap}_{F,p}(K) = \inf \{ \int_{\R^n} F^p(Dv) dx: v \in C_c^{\infty}(\R^n), v \geq 1\text{ on }K \},
\end{align*}
where $C_c^\infty(\R^n)$ is the set of smooth functions with compact support in $\R^n$.

For completeness, let us also introduce the anisotropic $p$-capacity for a general set in $\R^n$ as follows. For an open set $U\subset \R^n$, define
\begin{align*}
\mathrm{Cap}_{F,p}(U)=\sup_{K\subset U\;\text{compact}}\mathrm{Cap}_{F,p}(K).
\end{align*}
Then for any set $E\subset \R^n$, define
\begin{align*}
\mathrm{Cap}_{F,p}(E)=\inf_{E\subset U\;\text{open}}\mathrm{Cap}_{F,p}(U).
\end{align*}
We remark that in this paper we are mainly concerned with the anisotropic $p$-capacity for a compact set $K$.

Next we recall some basic facts on the anisotropic $p$-capacity.
\begin{prop}[\cite{HKM06,Maz11}]\label{prop-property}
The set function $E\mapsto \mathrm{Cap}_{F,p}(E)$ for $E\subset \R^n$ enjoys the following properties.
\begin{enumerate}
  \item If $E_1\subset E_2$, then $\mathrm{Cap}_{F,p}(E_1)\leq \mathrm{Cap}_{F,p}(E_2)$.
  \item For a compact set $K$, we have $\mathrm{Cap}_{F,p}(K)=\mathrm{Cap}_{F,p}(\pt K)$.
  \item If $\{K_i\}_{i\geq 1}$ is a decreasing sequence of compact sets in $\R^n$ with $K=\cap_{i\geq 1}K_i$, then
  \begin{equation*}
  \mathrm{Cap}_{F,p}(K)=\lim_{i\rightarrow \infty}\mathrm{Cap}_{F,p}(K_i).
  \end{equation*}
\end{enumerate}
\end{prop}

Under some regularity assumptions on a compact set $K$, there exists a unique weak solution to the following partial differential equation
\begin{equation}\label{eq-capacitary-potential}
		\begin{cases}
			\Delta_{F, p} u = 0 &\text{ in }\R^n\setminus K, \\
			u = 1 &\text{ on }  K,\\
			u(x)\rightarrow  0&\text{ as } |x| \to \infty.
		\end{cases}
	\end{equation}
The weak solution $u$ is called the \emph{anisotropic $p$-capacitary potential} of $K$. Here a function $u\in W_{loc}^{1,p}(U)$ for an open set $U\subset \R^n$ is called a \emph{weak solution} of $\Delta_{F, p} u=f$ on $U$ for $f\in L^q_{loc}(U)$ with $q=p/(p-1)$, if
\begin{equation*}
-\int_U\langle F^{p-1}DF(Du),Dv\rangle dx=\int_U fvdx
\end{equation*}
for any $v\in C_c^\infty(U)$. In the literature $\Delta_{F, p} u$ is called the \emph{anisotropic $p$-Laplacian} of $u\in W_{loc}^{1,p}(U)$. For a $C^2$ function $u$, at its regular points (where $Du\neq 0$) we have
\begin{equation*}
\Delta_{F, p} u :=\frac{1}{p}\mathrm{div}(D (F^p)(D u))= F^{p-2}( FF_{ij} + (p-1) F_i F_j) u_{ij}.
\end{equation*}
When $p=2$, the operator $\Delta_{F} u:=\Delta_{F, 2} u$ is called the \emph{anisotropic Laplacian} of $u\in W_{loc}^{1,2}(U)$.

For later use, we collect some results from \cite{AGHLV17} on the anisotropic $p$-capacitary potential $u$ of a compact convex set.
\begin{prop}[Lemmas~4.1, 4.3 and 4.4 in \cite{AGHLV17}]\label{prop-capacitary-potential}
Let $K\subset \R^n$ be a compact convex set with $\mathrm{Cap}_{F,p}(K)>0$. Then there exists a unique locally H\"{o}lder continuous function $0<u\leq 1$ on $\R^n$ satisfying \eqref{eq-capacitary-potential} with $u\in L^{np/(n-p)}(\R^n)$ and $|Du|\in L^p(\R^n)$. For each $t\in (0,1)$, the set $\{x\in \R^n:u(x)>t\}$ is convex. Moreover, if the interior of $K$ is non-empty, then $Du(x)\neq 0$ for $x\in\R^n\setminus K$.
\end{prop}
\begin{rem}
When $K$ is of non-empty interior, we see $\mathrm{Cap}_{F,p}(K)>0$ in view of (1) in Proposition~\ref{prop-property}, since we can find a small closed ball in $K$ with positive anisotropic $p$-capacity.
\end{rem}

\subsection{Anisotropic geometry of hypersurfaces}\label{sec2.2}

In this subsection we review the anisotropic geometry of hypersurfaces in the Euclidean space which is classical in the differential geometry. In contrast, in Section~\ref{sec2.3} we will introduce the relatively new anisotropic geometry of Euclidean hypersurfaces in the geometric analysis, especially where an anisotropic curvature flow is considered.

Let $F$ be a Minkowski norm on $\R^n$. We can define its dual norm $F^0$ as follows.
\begin{defn}
The \emph{dual norm} $F^0$ of $F$ is defined as
\begin{equation*}
F^0(x) = \sup_{\xi\neq 0} \f{\langle \xi, x \rangle}{F(\xi)}.
\end{equation*}
\end{defn}
It is known that $F^0$ is also a Minkowski norm.

Recall that $F$ and $F^0$ satisfy the following properties, which are very useful when we want to understand the relationship between the unit normal $\nu$ and the anisotropic unit normal $\nu_F$ of a hypersurface below.
\begin{prop}
\begin{enumerate}
  \item $F(DF^0(x)) = 1$, $F^0(DF(\xi)) = 1$.
  \item $F^0(x) DF(DF^0(x)) = x$, $F(\xi) DF^0(DF(\xi)) = \xi$.
\end{enumerate}
\end{prop}

%Besides, we note the following result.
%\begin{prop}
%The following statements are equivalent:
%\begin{enumerate}
%  \item $F^2$ is strongly convex in $\mathbb{R}^n$;
%  \item $\left(\mathrm{Hess}_{\mathbb{S}^{n-1}} F + F g_{\mathbb{S}^{n-1}}\right)|_{p}$ is positive definite on $p \in \mathbb{S}^{n-1}$.
%\end{enumerate}
%\end{prop}

Next we define the Wulff ball and the Wulff shape determined by $F$.
\begin{defn}
The \emph{Wulff ball} $\mathcal{W}$ centered at the origin is defined as
\begin{equation*}
\mathcal{W} := \{ x\in\mathbb{R}^{n}: F^0(x) < 1 \}.
\end{equation*}
Its boundary $\partial \mathcal{W}$ is called the \emph{Wulff shape}.
\end{defn}

%\begin{rem}
%	use 1-homogeneous properties of $F$, if $x \in \mathbb{S}^n$,   ontained $DF(x) = \nabla^{\mathbb{S}^n} F + \frac{\partial F}{\partial x} x = \nabla^{\mathbb{S}^n} F + Fx$.
%\end{rem}

Given a Wulff ball $\mathcal{W}$, we can recover $F$ as the support function of $\mathcal{W}$, namely,
\begin{align*}
F(\xi) = \sup_{X \in \mathcal{W}} \langle \xi, X \rangle,\quad  \xi\in \mathbb{S}^{n-1}.
\end{align*}

Now we introduce the anisotropic area of a smooth oriented hypersurface $X:N\to M\subset \R^n$.
\begin{defn}
	Let $M \subset \mathbb{R}^n$ be a smooth oriented hypersurface and $\nu$ be its unit normal vector. We define the \emph{anisotropic area} of $M$ as $|M|_F := \int_{M} F(\nu) \,d\mu$. Denote by $d\mu_F=F(\nu)d\mu$ the anisotropic area element of $M$.
\end{defn}

\begin{rem}
For $M=\pt \mathcal{W}$, we can check by the divergence theorem that
\begin{equation*}
|\partial \mathcal{W}|_F := \int_{\partial \mathcal{W}} F(\nu) \,d\mu = \int_{\partial \mathcal{W}} \langle X, \nu\rangle \,d\mu = \int_{\mathcal{W}} \mathrm{div} X\,dx = n|\mathcal{W}|.
\end{equation*}
\end{rem}

Next we introduce the anisotropic Gauss map for an oriented hypersurface in $\R^n$.
\begin{defn}
The \emph{anisotropic Gauss map} $\nu_F: M \to \pt \mathcal{W}$ from an oriented hypersurface $M$ in $\R^{n}$ to the Wulff shape $\pt\mathcal{W}$ is defined by
\begin{align*}
\nu_F : \ & M \to \pt \mathcal{W},\\
&X \mapsto DF(\nu(X))= F(\nu(X)) \nu(X) + \nabla^{\mathbb{S}^{n-1}} F (\nu(X)),
\end{align*}
where $\nu$ is the unit normal vector of $M$.
\end{defn}
\begin{rem}
The vector $\nu_F$ is also called the anisotropic unit normal of the hypersurface.
\end{rem}

Let $A_F$ be the 2-tensor on $\SS^{n-1}$ defined by
\begin{align*}
A_F(\xi) =( \nabla^{\mathbb{S}^{n-1}})^2 F(\xi) + F(\xi) \sigma,\quad \xi \in \SS^{n-1},
\end{align*}
where $\sigma$ is the standard metric on $\mathbb{S}^{n-1}$.

\begin{defn}
The \emph{anisotropic principal curvatures} $\kappa_1^F, \dots, \kappa_{n-1}^F$ of a smooth oriented hypersurface $M$ in $\R^n$ are defined as the eigenvalues of the tangent map
\begin{equation*}
d\nu_F:T_XM\to T_{\nu_F(X)}\PW\cong T_XM.
\end{equation*}
The \emph{anisotropic mean curvature} is defined as
\begin{align*}
H_F := \mathrm{tr}(d\nu_F) =\sum_i \kappa_i^F = \sum_{i, j, k} (A_F)_{i}^{j} \left(\nu(X)\right) g^{ik}(X) h_{kj}(X),
\end{align*}
where $g$ and $h$ are the first and second fundamental forms of the hypersurface respectively. The \emph{anisotropic Gaussian curvature} is defined as
\begin{align*}
K_F := \mathrm{det}(d\nu_F) = \prod_i \kappa_i^F=\det(A_F)\det(g^{-1}h).
\end{align*}
Given an integer $1\leq k\leq n-1$, the \emph{$k$th anisotropic mean curvature} is defined as
\begin{align*}
\sigma_k(\kappa^F) := \sum_{1\le i_1 < \cdots < i_k \le n-1 } \kappa_{i_1}^F \cdots \kappa_{i_k}^F.
\end{align*}
\end{defn}
Note that $H_F=\sigma_1(\kappa^F)$ and $K_F=\sigma_{n-1}(\kappa^F)$.

\subsection{Anisotropic Riemannian metric on hypersurfaces and inverse anisotropic mean curvature flow}\label{sec2.3}

The materials presented in Section~\ref{sec2.2} mostly suffice for anisotropic problems in the differential geometry. However, when we come to the field of the geometric analysis, especially when we consider anisotropic curvature flows, those seem insufficient.

In \cite{And01}, to study the volume-preserving anisotropic curvature flow in the Euclidean space, Ben~Andrews introduced a new Riemannian metric on hypersurfaces in $\R^n$ associated with the anisotropy $F$. Later Chao Xia \cite{Xia13} reformulated Andrews' setting and used it to study the inverse anisotropic mean curvature flow \cite{Xia17}. It turns out that Andrews' new metric is very useful in the geometric analysis concerning the anisotropic geometry of Euclidean hypersurfaces. In this subsection, we shall review the settings in \cite{And01} and \cite{Xia13} and also the inverse anisotropic mean curvature flow in \cite{Xia17}.

Using the dual norm $F^0$, for any non-vanishing vector field $z(x)$ on $\R^n$ we can define a Riemannian metric $G$ on $T\mathbb{R}^n$ by
\begin{align*}
 	G(z)  &(V, W):= \sum_{\alpha, \beta = 1}^{n} \frac{\partial^2 ((F^0)^2(z)/2)}{\partial z^{\alpha} \pt z^{\beta}} V^{\alpha} W^{\beta}, \\
 & \quad \forall 0\neq z(x) \in \R^{n},\quad  V, W \in T_{x} \R^{n}.
 \end{align*}

In general, the third derivatives of $(F^0)^2(z)$ do not vanish. So we can define a $3$-tensor $Q$ as
\begin{align*}
	Q(z) ( U, &V, W) := \sum_{\alpha, \beta, \gamma = 1}^{n} \frac{\partial^3 ((F^0)^2 (z)/2) }{\partial z^\alpha \partial z^{\beta} \partial z^{\gamma}} U^{\alpha} V^{\beta} W^{\gamma},\\
& \forall 0\neq z(x) \in \R^{n},\quad U, V, W \in T_{x} \R^{n}.
\end{align*}

By virtue of the $1$-homogeneity of $F^0$, we can derive
\begin{align*}
G(z) (z, z) = 1, \quad  G(z) (z, V) = 0 \quad\text{for }z \in \pt\mathcal{W},\quad V \in T_{z}\pt \mathcal{W},
\end{align*}
and
\begin{align*}
Q(z) (z, V, W) = 0 \quad \text{for } z \in \pt\mathcal{W},\quad V, W \in T_{z}\R^n.
\end{align*}
Let $M=X(N)$ be a smooth oriented hypersurface $X:N\to \R^n$ from a manifold $N$. Setting $z(X)=\nu_F(X)$, we get along the hypersurface
\begin{align*}
&\qquad G(\nu_F(X)) (\nu_F(X),\nu_F(X)) = 1,\\
&G(\nu_F(X)) (\nu_F(X), V) = 0 \quad\text{for } V \in T_XM.
\end{align*}
Motivated by the above identities, we define a Riemannian metric $\hat{g}$ on the hypersurface $M$ by restricting $G(\nu_F(X))$ on its tangent space,
\begin{align*}
\hat{g}(X) = G(\nu_F(X))|_{T_X M}, \quad X \in M.
\end{align*}

Let $\hat{\nabla}$ be the Levi-Civita connection of $\hat{g}$ on $M$. Then we have the following expressions for the anisotropic Riemannian metric $\hat{g}$, the anisotropic second fundamental form $\hat{h}$ and the $3$-tensor $Q$:
\begin{align*}
&\;\; \hat{g}_{ij} := G(\nu_F(X))(\partial_i X, \partial_j X),\\
&\hat{h}_{ij} := -G(\nu_F(X))\left(\nu_F, \pt_i\partial_j X\right),\\
&\, Q_{ijk} := Q(\nu_F)(\partial_i X, \partial_j X, \partial_k X),
\end{align*}
where $\{\pt_i\}_{i=1}^{n-1}$ are local coordinate vectors on $N$. For more details on $\hat{g}$, $\hat{h}$ and $Q$, see \cite{And01} and \cite{Xia17}.

In \cite{Xia17}, following the (isotropic) works \cite{Ger90,Urb90}, Chao Xia considered the inverse anisotropic mean curvature flow for a star-shaped $F$-mean convex hypersurface, i.e., the parabolic evolution equation
\begin{align}\label{eq-IAMCF}
\partial_t X = \frac{1}{H_F} \nu_F,
\end{align}
where $X:N\times [0,T)\to \R^n$ is a smooth family of embeddings from a closed manifold $N$ to $\R^n$. He proved the following result on this flow.
\begin{thm}[\cite{Xia17}]
Let $M$ be a smooth compact star-shaped and $F$-mean convex hypersurface without boundary in $\R^n$ ($n\geq  2$). Then the inverse anisotropic mean curvature flow starting from $M$ exists for all time, and converges smoothly and exponentially to an expanded Wulff shape determined by the initial hypersurface $M$.
\end{thm}
In this paper we also need the following evolution equations along the flow~\eqref{eq-IAMCF}.
\begin{prop}[Proposition 4.1 in \cite{Xia17}]
Along the flow \eqref{eq-IAMCF} we have
\begin{align}
%&\qquad \qquad \qquad \qquad\partial_t \nu_F = -\hat{\nabla}\left(\frac{1}{H_F}\right),\\
&\qquad \qquad \qquad \qquad \quad \partial_t d\mu_F = d\mu_F,\label{eq-area-element}\\
&\partial_t H_F = \frac{1}{H_F^2} \left( \hat{\Delta}H_F + \hat{g}^{ik} A_{pik} \hat{\nabla}^p H_F \right) - 2 \frac{|\hat{\nabla} H_F|_{\hat{g}}^2}{H_F^3} -\frac{1}{H_F}|\hat{h}|_{\hat{g}}^{2},\label{eq-mean-curvature}
\end{align}
where $\hat{\Delta}$ is the Laplacian corresponding to $\hat{g}$ and
\begin{equation*}
A_{ijk} := -   (\hat{h}_i^l Q_{jkl} + \hat{h}_j^l Q_{ikl} - \hat{h}_k^l Q_{ijl} )/2.
\end{equation*}
\end{prop}

Last we recall the formula for the integration by parts in the anisotropic setting.
\begin{prop}[{\cite[Lemma~2.8]{Xia13}}, {\cite[Lemma~3.2]{Xia17}}]
Let $f_1$ and $f_2$ be two smooth functions on a compact hypersurface $M\subset \R^n$ without boundary. Then we have
	\begin{equation}\label{eq-integration}
		\int_{M} f_1\left(\hat{\Delta} f_2 + \hat{g}^{jk} A_{ijk} \hat{\nabla}^i f_2\right)\,d\mu_F = -\int_{M} \hat{\nabla}^i f_1 \hat{\nabla}_i f_2\,d\mu_F.
	\end{equation}
\end{prop}

\section{Proof  of Theorem \ref{thm1}}\label{sec3}

In this section we first introduce a general approach which will be employed in the proofs of Theorems~\ref{thm1}, \ref{thm2} and \ref{thm3}, and then we present the proof of Theorem~\ref{thm1}.

\subsection{General approach}
	The general approach presented here may date back to \cite{PS51,Sze31} and has been applied successfully in some works, e.g., in \cite{FS14,Xiao16,Xiao17,BM08}.
	
	Let $U=\R^n\setminus K$. For the proofs of Theorems~\ref{thm1}, \ref{thm2} and \ref{thm3}, we use the fact that $\overline{U}$ can be foliated by a family of hypersurfaces $M_t$ ($t\geq 0$) such that $M_0=\pt U=\pt K$ and $M_t\rightarrow \infty$ as $t\rightarrow \infty$. Later we shall specify which foliation we choose in each case. For the moment we note that $M_t$ defines a function $\psi$ on $\overline{U}$ with level sets being $M_t$, i.e.,
	\begin{equation*}
		\psi(x)=t, \text{ when } x\in M_t.
	\end{equation*}
	Define two functions
	\begin{align*}
		&T_{p}(t):=\int_{\psi(x)=t}\frac{F^p(D\psi)}{|D \psi|} d\mu,\\
		&\lambda(t):=\frac{\int_t^\infty  T_{p}^{-1/(p-1)}(s)ds}{\int_0^\infty  T_{p}^{-1/(p-1)}(s)ds}.
	\end{align*}
	Note that $\lambda\in C^1([0,\infty))$ is a non-increasing function satisfying $\lambda(0)=1$ and $\lambda(\infty)=0$. Moreover, set
	\begin{equation*}
		f(x)=\lambda(\psi(x)),\quad x\in \overline{U}.
	\end{equation*}

Now we point out that in the definition of $\mathrm{Cap}_{F,p}(K)$, we can actually choose test functions as Lipschitz or even locally Lipschitz functions $f(x)$ satisfying $f=1$ on $\pt K$ and $f(x)\to 0$ as $|x|\to \infty$. See e.g. \cite[Definition~1]{BM08} and the remark below \cite[Definition~1.2]{Shu}.

	Back to our setting we can check in each case that the function $f$ we choose later is a Lipschitz function satisfying $f=1$ on $\pt K$ and $f(x)\to 0$ as $|x|\to \infty$. So $f$ is admissible for the definition of $\mathrm{Cap}_{F,p}(K)$ and we get
	\begin{equation*}
		\mathrm{Cap}_{F, p}(K)\leq \int_{U}F^p(D f)  dx.
	\end{equation*}
	Next we estimate $\int_{U}F^p(D f) dx$. By the co-area formula we obtain
	\begin{align*}
		\int_{U}F^p(Df) dx&=\int_{U}|\lambda'(\psi(x))|^p F^p(D\psi) dx\\
		&=\int_0^\infty |\lambda'(t)|^p \int_{\psi(x)=t}\frac{F^p(D\psi)}{|D \psi|} d\mu\, dt\\
		&=\int_0^\infty |\lambda'(t)|^p T_{p}(t)dt.
	\end{align*}
	On the other hand, by the H\"{o}lder inequality we have
	\begin{align*}
		1&=(\lambda(0))^p=\left(-\int_0^\infty \lambda'(t)dt\right)^p\\
		&=\left(\int_0^\infty (-\lambda'(t))T_{p}^{1/p}(t)\cdot T_{p}^{-1/p}(t)dt\right)^p\\
		&\leq \int_0^\infty |\lambda'(t)|^p T_{p}(t)dt  \cdot \left(\int_0^\infty  T_{p}^{- \frac{1}{p-1}}(t)dt\right)^{p-1},
	\end{align*}
	with the equality if $\lambda'(t)=cT_{p}^{-1/(p-1)}(t)$. In view of the definition of $\lambda$, the above inequality is indeed an equality. As a consequence, we obtain
	\begin{equation*}
		\int_{U}F^p(Df) dx= \left(\int_0^\infty  T_{p}^{-\frac{1}{p-1}}(t)dt\right)^{1-p},
	\end{equation*}
	and then
	\begin{equation}\label{eq-general}
		\mathrm{Cap}_{F, p}(K)\leq \left(\int_0^\infty  T_{p}^{-\frac{1}{p-1}}(t)dt\right)^{1-p} .
	\end{equation}

\subsection{Proof of Theorem~\ref{thm1}}

Now we consider the inverse anisotropic mean curvature flow
\begin{equation*}
\begin{cases}
X:N\times [0,T)\to \R^n,\\
\partial_t  X = \dfrac{1}{H_F}\nu_F,
\end{cases}
\end{equation*}
with the initial condition $X(N, 0) = \partial K$. % we write level set function of flow by $u(x): \R^n \to \R$ s.t. $u(X(p, t)) = t$ for any $t > 0$. write $\psi: \R \to \R$ be a function to be choosen, and satisfies $\psi(0) = 1, \psi(\infty) = 0$.

%Clearly, for any $f$ satisfies $f|_{\partial K} = 1, f(\infty) = 0$

%\[
%Cap_{F, p}(K) \le \int_{\R^n \setminus K} F(Df)^p \,dV
%\]
%
%we write $f = \psi \circ u$, then RHS of above is
%
%\begin{equation*}
%	\begin{aligned}
%		Cap_{F, p}(K) &\le \int_{0}^{\infty} |\psi'(t)|^p \int_{\{x| u(x) = t\}}\frac{F(Du)^p}{|Du|}\,dS_t\,dt \\
%	\end{aligned}
%\end{equation*}
%
%Define $U_p(t) = \int_{\{x| u(x) = t\}}\frac{F(Du)^p}{|Du|}\,dS_t$, and $\psi(t) = \frac{\int_{t}^{\infty} U_p(s)^{\frac{1}{1-p}}ds}{\int_{0}^{\infty} U_p(s)^{\frac{1}{1-p}}ds}$. then take in above inequality, we obtained
%
%\begin{equation*}
%	\begin{aligned}
%		Cap_{F, p}(K) \le \left( \int_{0}^{\infty} U_p(t)^{\frac{1}{1-p}}\,dt \right)^{1-p}
%	\end{aligned}
%\end{equation*}

In this case we have $\psi(X(p,t))=t$. Taking the derivative with respect to $t$ yields
\begin{equation*}
\langle D\psi,\frac{1}{H_F}\nu_F\rangle =1.
\end{equation*}
Noting $\nu_F=F(\nu)\nu+\nabla^{\SS^{n-1}}F(\nu)$ and $D\psi=|D\psi|\nu$, we obtain
\begin{align*}
D\psi = \frac{H_F}{F(\nu)} \nu.
\end{align*}
So we get
\begin{align*}
T_p(t) &= \int_{ \psi(x) = t} \frac{F^p(D\psi)}{|D\psi|}\,d\mu= \int_{  \psi(x) = t} H_F^{p-1} F(\nu) \,d\mu\\
&=\int_{  \psi(x) = t} H_F^{p-1} \,d\mu_F.
\end{align*}
Using the evolution equations \eqref{eq-area-element} and \eqref{eq-mean-curvature} along the inverse anisotropic mean curvature flow, we derive
\begin{align}
	&\frac{d}{dt} T_p(t)	= \int_{ \psi(x) = t } \bigg((p-1) H_F^{p-2}\bigg( \frac{1}{H_F^2}(\hat{\Delta}H_F + \hat{g}^{ik} A_{pik} \hat{\nabla}^p H_F)\nonumber\\
 &\qquad -2 \frac{|\hat{\nabla}H_F|^2_{\hat{g}}}{H_F^3} - \frac{|\hat{h}|^2_{\hat{g}}}{H_F} \bigg) + H^{p-1}_F \bigg) \,d\mu_F\nonumber\\
 &=\int_{  \psi(x) = t } \bigg(-(p-1)(p-2) H_F^{p-5}|\hat{\nabla}H_F|^2_{\hat{g}} - (p-1)H_F^{p-3}|\hat{h}|^2_{\hat{g}}  + H^{p-1}_F \bigg) \,d\mu_F\nonumber\\
 &\leq \int_{  \psi(x) = t } \bigg(-(p-1)(p-2) H_F^{p-5}|\hat{\nabla}H_F|^2_{\hat{g}} + \frac{n-p}{n-1}H_F^{p-1}  \bigg) \,d\mu_F, \label{eq-computation}
 \end{align}
where in the second equality we used \eqref{eq-integration} for the integration by parts and in the inequality we used $|\hat{h}|^2_{\hat{g}} \geq  H_F^2/(n-1)$.

First assume $p\in [2,n)$. We get
\begin{align}\label{eq-inequality}
	\frac{d}{dt} T_p(t)	\leq &\int_{\psi(x) = t} \frac{n-p}{n-1} H_F^{p-1}\,d\mu_F =\frac{n-p}{n-1}T_p(t).
\end{align}
Solving the above ordinary differential inequality we get
\begin{equation}\label{eq-Up}
	T_p(t) \le e^{\frac{n-p}{n-1} t} \int_{\partial K}H_F^{p-1}\,d\mu_F,\quad \forall p\in [2,n).
\end{equation}

Now we consider the case $2\leq p< n$ in Theorem~\ref{thm1}. Plugging \eqref{eq-Up} into Inequality~\eqref{eq-general} leads to
\begin{align*}
	\mathrm{Cap}_{F, p}(K) \le \left(\frac{(p-1)(n-1)}{n-p}\right)^{1-p} \int_{\partial K} H_F^{p-1}\,d\mu_F.
\end{align*}

Next we consider the case $1<p\leq 2\leq q<n$. We obtain by the H\"{o}lder inequality
\begin{align*}
	T_p(t)& = \int_{\psi(x) = t} H_F^{p-1}\,d\mu_F \\
&\le \left( \int_{\psi(x) = t} H_F^{q-1}\,d\mu_F \right)^{\frac{p-1}{q-1}} |M_t|_F^{\frac{q-p}{q-1}}\\
&=T_q^{\frac{p-1}{q-1}}(t) |M_t|_F^{\frac{q-p}{q-1}}.
\end{align*}
Noticing $|M_t|_F =  |\partial K|_F \, e^t$ and using \eqref{eq-Up} for $T_q(t)$, we get
\begin{align*}
		T_p(t) &\le \left(\int_{\partial K} H_F^{q-1}\,d\mu_F\right)^{\frac{p-1}{q-1}} |\partial K|_F^{\frac{q-p}{q-1}} e^{\frac{(n-q)(p-1)}{(n-1)(q-1)}t + \frac{q-p}{q-1}t}\\
		&= \left(\int_{\partial K} H_F^{q-1}\,d\mu_F\right)^{\frac{p-1}{q-1}} |\partial K|_F^{\frac{q-p}{q-1}} e^{\frac{n-p}{n-1}t}.
\end{align*}
Therefore in light of Inequality~\eqref{eq-general} we conclude
\begin{align*}
		\mathrm{Cap}_{F, p}(K) &\le \left(\int_{\partial K} H_F^{q-1}\,d\mu_F\right)^{\frac{p-1}{q-1}} |\partial K|_F^{\frac{q-p}{q-1}}\left(\int_{0}^{\infty} e^{\frac{(n-p) t}{(n-1)(1-p)}}\,dt\right)^{1-p} \\
		&= \left( \frac{(p-1)(n-1)}{n-p} \right)^{1-p} \left(\int_{\partial K} H_F^{q-1}\,d\mu_F\right)^{\frac{p-1}{q-1}} |\partial K|_F^{\frac{q-p}{q-1}}.
\end{align*}

If the equality holds, then \eqref{eq-inequality} must be an equality. So $M_t$ is anisotropicly umbilical. By the anisotropic Codazzi equation (see e.g. \cite{And01,Xia13,Xia17}), $H_F$ is constant, which implies that $M_t$ is a translated scaled Wulff shape (\cite{HLMG09}). Hence $M$ is a translated scaled Wulff shape.

\section{Proof of Theorem \ref{thm2}}\label{sec4}

Given a compact oriented surface $\Sigma$ without boundary in $\R^3$, we define the \emph{anisotropic Hawking mass} $m_H^F(\Sigma)$ of $\Sigma$ by
\begin{align*}
	m_H^F(\Sigma) := \sqrt{\frac{|\Sigma|_F}{4 |\PW|_F}} \left(1 - \frac{\int_{\Sigma} H_F^2\,d\mu_F}{4|\PW|_F}\right),
\end{align*}
where as before $|\cdot|_F$ denotes the anisotropic area, $H_F$ is the anisotropic mean curvature, and $d\mu_F = F(\nu)d\mu$ is the anisotropic area element. As mentioned in the Introduction, the classical (isotropic) Hawking mass $m_H(\Sigma)$ (corresponding to $F(\xi)=|\xi|$) was introduced by S.~W.~Hawking \cite{Haw68}.

First we have the following observation.
\begin{lem}\label{Gauss-curvature}
For any smooth compact star-shaped hypersurface $M$ without boundary in $\R^n$ with $K_F$ denoting its anisotropic Gaussian curvature, we have
\begin{align*}
\int_{M} K_F \,d\mu_F = |\PW|_F.
\end{align*}
\end{lem}

\begin{proof}
	%Since the integration of anisotropic Gauss Curvature is invariant under any variation, see\cite{reilly1976}. take the variation be $\frac{1}{H_F}$, \cite{Xia17} shows $\Sigma$ evolves to a strictly convex hypersurface $\Sigma'$
	
Since both $M$ and the Wulff shape $\pt \mathcal{W}$ are star-shaped, they can be represented by using the polar coordinates as
\begin{align*}
M&=\{(r_1(p),p):p\in \SS^{n-1}\},\\
\pt \mathcal{W}&=\{(r_2(p),p):p\in \SS^{n-1}\},
\end{align*}
where $r_i\in C^\infty(\SS^{n-1})$ ($i=1,2$) are two smooth positive functions on $\SS^{n-1}$.

Now consider the smooth variation $X:\SS^{n-1}\times [0,1]$ connecting $M$ and $\pt \mathcal{W}$,
\begin{align*}
X(p,t)=((1-t)r_1(p)+tr_2(p),p)\in \R^n,\quad p\in\SS^{n-1},\quad t\in [0,1],
\end{align*}
and let $M_t=X(\SS^{n-1},t)$ for $t\in [0,1]$.

Since the anisotropic integration of the anisotropic Gaussian curvature over a compact smooth hypersurface without boundary is invariant under any smooth variation (\cite[Theorem~4]{Rei76}), we get
\begin{align*}
\int_{M} K_F \,d\mu_F  =\int_{M_t} K_F \,d\mu_F= \int_{\pt \mathcal{W}}K_F \,d\mu_F &= |\PW|_F.
\end{align*}
\end{proof}

The next result we will use later is a sharp lower bound for the anisotropic Willmore energy of an embedded compact hypersurface without boundary. It may be well-known to specialists. See e.g. \cite[Theorem~1.2]{XY20} as a corollary of the main result there. Here we present a direct proof for the readers' convenience.
\begin{prop}\label{willmoreineq}
Let $M$ be an embedded smooth compact oriented hypersurface without boundary in $\R^n$. We have
\begin{align*}
\int_{M} |H_F|^{n-1} d\mu_F\ge (n-1)^{n-1} |\PW|_F,
\end{align*}
with the equality holding if and only if $M$ is a translated scaled Wulff shape.
\end{prop}

\begin{proof}
Fix any $p\in \pt\mathcal{W}$. Consider the tangent hyperplane $T_p\pt \mathcal{W}$. By moving this hyperplane from the infinity to $M$, we can find a point $x\in M$ such that $T_xM$ is parallel to $T_p\pt \mathcal{W}$ and a neighbourhood of $x$ on $M$ lies on the same side of $T_xM$ from which $\nu_F(x)=p$ points. Note that at $x$, all the anisotropic principal curvatures are nonnegative. Let $\widehat{M}$ be the set of points on $M$ where all the anisotropic principal curvatures are nonnegative. Then in view of the above observation and the Sard's theorem, we can readily derive
\begin{align*}
\int_{M} |H_F|^{n-1}d\mu_F &\ge \int_{\widehat{M}} |H_F|^{n-1}d\mu_F \\
&\ge (n-1)^{n-1}\int_{\widehat{M}} K_F d\mu_F\\
&\ge (n-1)^{n-1} |\PW|_F.
\end{align*}

Next assume that the equality holds. By applying the inequality to each connected component of $M$, we know that $M$ itself must be connected. Then analyzing the equality cases of the above sequence of inequalities, we see that all points on $M$ are anisotropicly umbilical. Then the same argument as at the end of Section~\ref{sec3} implies that $M$ is a translated scaled Wulff shape. The proof is complete.
\end{proof}

From the above result, we have an immediate corollary.
\begin{corr}\label{corr-non-positive}
Let $\Sigma$ be an embedded smooth oriented compact surface without boundary in $\R^{3}$. Then $m_{H}^F(\Sigma)\leq 0$, with the equality holding if and only if $\Sigma$ is a translated scaled Wulff shape.
\end{corr}

%\begin{proof}
%	Easily from proposition \ref{willmoreineq} if $n = 2$.
%\end{proof}

Next we prove a monotonicity result for the anisotropic Hawking mass of a star-shaped $F$-mean convex surface in $\R^3$ along the inverse anisotropic mean curvature flow. The corresponding and more general result in the isotropic case can be found in \cite{Gero73,HI01}.
\begin{prop}
Let $\Sigma$ be a compact star-shaped $F$-mean convex surface without boundary in $\R^3$. Along the inverse anisotropic mean curvature flow \eqref{eq-IAMCF} starting from $\Sigma$, the anisotropic Hawking mass $m_H^F(\Sigma_t)$ is non-decreasing in $t$. Moreover, if $(d/dt)m_H^F(\Sigma_t)=0$ at some time $t>0$, then $\Sigma$ is a translated scaled Wulff shape.
\end{prop}

\begin{proof}
First recalling the computation in \eqref{eq-computation}, we get (let $p=3$ there)
	\begin{align*}
					&\frac{d}{d t} \int_{\Sigma_t} H_F^2\,d\mu_F= \int_{\Sigma_t}\left( -2 \frac{|\hat{\nabla} H_F|_{\hat{g}}^2}{H_F^2} -2|\hat{h}|_{\hat{g}}^{2} +H_F^2\right) d\mu_F.
	\end{align*}%
%where we used the integration by parts
%\begin{align*}
%			\int_{\Sigma_t} \frac{2}{H_F} \left( \hat{\Delta}H_F + \hat{g}^{ik} A_{pik} \hat{\nabla}_p H_F \right)\,dS_F &= -\int_{\Sigma_t} \hat{\nabla}\left(\frac{2}{H_F}\right) \cdot \hat{\nabla}H_F\,dS_F\\
%			&=\int_{\Sigma_t} 2\frac{|\hat{\nabla}H_F|_{\hat{g}}^2}{H_F^2}\,dS_F.
%\end{align*}
Using $H_F^2 =2K_F+ |\hat{h}|_{\hat{g}}^2$ and Lemma~\ref{Gauss-curvature}, we obtain
\begin{align*}
		\frac{d}{d t} \int_{\Sigma_t} H_F^2\,d\mu_F &= \int_{\Sigma_t} \left(-2 \frac{|\hat{\nabla} H_F|_{\hat{g}}^2}{H_F^2} -|\hat{h}|_{\hat{g}}^{2} + 2K_F \right)d\mu_F\\
		&\le 2|\PW|_F - \int_{\Sigma_t} |\hat{h}|_{\hat{g}}^{2} \,d\mu_F.
\end{align*}
Noting $|\hat{h}|_{\hat{g}}^{2} \geq H_F^2/2$, we get
\begin{equation*}
		\frac{d}{d t} \int_{\Sigma_t} H_F^2\,d\mu_F \le \frac{1}{2} \left(4|\PW|_F - \int_{\Sigma_t} H_F^2\,d\mu_F\right).
\end{equation*}
Therefore in view of $(d/dt)|\Sigma_t|_F=|\Sigma_t|_F$ we conclude
\begin{align*}
	&\frac{d}{dt} \left(  |\Sigma_t|_F^{\frac{1}{2}} \left( 4|\PW|_F - \int_{\Sigma_t} H_F^2 \,d\mu_F \right) \right) \\
&\ge \frac{1}{2}|\Sigma_t|_F^{\frac{1}{2}} \left( 4|\PW|_F - \int_{\Sigma_t} H_F^2 \,d\mu_F \right)-\frac{1}{2}|\Sigma_t|_F^{\frac{1}{2}} \left( 4|\PW|_F - \int_{\Sigma_t} H_F^2 \,d\mu_F \right)\\
&=0.
\end{align*}
So the anisotropic Hawking mass is non-decreasing in $t$.

Now assume $(d/dt)m_H^F(\Sigma_t)=0$ at some time $t>0$. Then checking the above argument we see that $H_F$ is constant on $\Sigma_t$, which implies that $\Sigma_t$ is a translated scaled Wulff shape (\cite{HLMG09}). So the initial surface $\Sigma$ is a translated scaled Wulff shape. The proof is complete.

\end{proof}

%
%	We give the definition of Anisotropic P-capacity. for a smooth bounded domain $K \subset \left(\R^3, F\right)$, we denoted Anisotropic P-capacity of $K$ by $Cap_{PF}\left(K\right)$, it has value
	
%	\begin{equation*}
%		Cap_{PF} \left(K\right) = \inf_u \{ \int_{\R^3\setminus K} F(Du)^p \,dV_{\R^3}| u \in C_{0}^{\infty} \left(K\right) , u|_{\partial K} = 1 \}
%	\end{equation*}

\begin{proof}[Proof of Theorem \ref{thm2}]
%	We also use anisotropic inverse mean curvature flow to foliation $\R^3 \setminus K$. denoted test function is $u = f \circ \psi$. ,
%	
%	
%	By definition $Cap_F(K) = \inf_v\int_{\R^3 \setminus K} F(Dv)^p \,dV_{\R^3}$, where $v$ satisfies $v|_{\Sigma} = 0, v|_{\infty} =1$.
	
	%Use inverse anisotropic mean curvature flow $\partial_t X = \frac{\nu_F}{H_F}$ to foliation $\R^3\setminus K$. define $\psi: \R^3\setminusK \to \R^{+}$ be time level set of inverse anisotropic mean curvature flow s.t. if $\forall x \in X(\cdot, t)$, then $\psi(x) = t$. and define $f: \R \to \R$ s.t. $f(0) = 0, f(\infty) =1$.
%	
%	Test function $v$ we choose by $v(x) = f(\psi(x))$, then
%	
%	\begin{equation*}
%		\int_{\R^3 \setminus K} F(Dv)^p \,dV_{\R^3} = \int_{0}^{\infty} \left(f'(t)\right)^p \int_{\{x| \psi(x) = t\}} \frac{F(D\psi)^p}{|D \psi|}\,dS_t\,dt
%	\end{equation*}
%	
%Then
%	
%	\begin{equation*}
%		Cap_{F, p}\left(K\right) \le \int_{0}^{\infty} |f'(t)|^p \int_{\{x| \psi(x) = t\}} \frac{F(D\psi)^p}{|D\psi|}\,dS_t \,dt
%	\end{equation*}
%	
%	
%	Level set function $\psi$ satisfies $\psi\left( X(p, t) \right) = t$, differential on $t$ gives $D\psi \cdot \frac{\nu_F}{H_F} = 1$, gives $D\psi = \frac{H_F}{F(\nu)}\nu$. where $U_p(t) := \int_{\{x| \psi(x) = t\}} \frac{F(D\psi)^p}{|D\psi|}\,dS_t = \int_{\{x| \psi(x) = t\}} H_F^{p-1} F(\nu) \,dS_t = \int_{\{x| \psi(x) = t\}} H_F^{p-1} \,dS_t^F$.

For the proof we still use the inverse anisotropic mean curvature flow. Recall along this flow we have proved
\begin{equation*}
		\mathrm{Cap}_{F, p}(K)\leq \left(\int_0^\infty  T_p^{\frac{1}{1-p}}(t)dt\right)^{1-p},
	\end{equation*}
where
\begin{align}\label{eq-inequality-surface}
T_p(t) &= \int_{\Sigma_t} H_F^{p-1} \,d\mu_F \leq \left(\int_{\Sigma_t} H_F^2\,d\mu_F\right)^{\frac{p-1}{2}} \left( |\Sigma_t|_F \right)^{\frac{3-p}{2}},
\end{align}
after using the H\"{o}lder inequality.

Note that $m_H^F(\Sigma_t)$ is monotone non-decreasing in $t$. So
\begin{align*}
m_H^F\left(\Sigma\right) \le m_H^F\left( \Sigma_t\right) = \sqrt{\frac{|\Sigma_t|_F}{4 |\PW|_F}} \left(1 - \frac{\int_{\Sigma_t} H_F^2\,d\mu_F}{4|\PW|_F}\right),
\end{align*}
which means
\begin{align*}
\int_{\Sigma_t} H_F^2 \,d\mu_F \le 4 |\PW|_F\left( 1- \sqrt{\frac{4 |\PW|_F}{|\Sigma_t|_F}} m_H^F\left(\Sigma\right)\right).
\end{align*}
Consequently, we get	
\begin{align*}
T_p(t) &\le \left( 4 |\PW|_F\right)^{\frac{p-1}{2}}\left(1 - \sqrt{\frac{4 |\PW|_F}{|\Sigma_t|_F}} m_H^F\left(\Sigma\right) \right)^{\frac{p-1}{2}}  \left( |\Sigma_t|_F \right)^{\frac{3-p}{2}}\\
			&= \left( 4 |\PW|_F\right)^{\frac{p-1}{2}}\left(1 - \sqrt{\frac{4 |\PW|_F}{|\Sigma|_F e^t}} m_H^F\left(\Sigma\right) \right)^{\frac{p-1}{2}}  \left( |\Sigma|_F e^t \right)^{\frac{3-p}{2}}\\
&=(4|\PW|_F)^{\frac{p-1}{2}}|\Sigma|_F^{\frac{3-p}{2}}\left( 1 - \sqrt{\frac{4 |\PW|_F}{|\Sigma|_F  }} m_H^F(\Sigma)e^{-t/2} \right)^{\frac{p-1}{2}}   e^{\frac{3-p}{2}t},
\end{align*}
where we used $|\Sigma_t|_F = |\Sigma|_F e^t$.
	
So we get
\begin{align*}
&\mathrm{Cap}_{F, p}(K)\leq \left(\int_0^\infty  T_p^{\frac{1}{1-p}}(t)dt\right)^{1-p}\\
&\leq (4|\PW|_F)^{\frac{p-1}{2}}|\Sigma|_F^{\frac{3-p}{2}}\left(\int_0^\infty \left( 1 - \sqrt{\frac{4 |\PW|_F}{|\Sigma|_F  }} m_H^F(\Sigma)e^{-t/2} \right)^{-1/2}e^{\frac{3-p}{2(1-p)}t}dt\right)^{1-p}.
\end{align*}	

If $m_H^F(\Sigma)=0$, then by Corollary~\ref{corr-non-positive} the surface $\Sigma$ is a translated scaled Wulff shape $r_0\pt \mathcal{W}+x_0$ ($r_0>0$, $x_0\in \R^3$) and we can compute directly to get
\begin{align*}
\mathrm{Cap}_{F, p}\left(K\right) = \left(\frac{3-p}{p-1}\right)^{p-1} |\PW|_F \; r_0^{3-p}.
\end{align*}
%
%\begin{align*}
%\mathrm{Cap}_{F, p}(K)&\leq  (4|\PW|_F)^{\frac{p-1}{2}}|\Sigma|_F^{\frac{3-p}{2}}\left(\int_0^\infty e^{\frac{3-p}{2(1-p)}t}dt\right)^{1-p}\\
%&=\left(\frac{3-p}{p-1}\right)^{p-1} |\PW|_F^{\frac{p-1}{2}} |\Sigma|_F^{\frac{3-p}{2}}.
%\end{align*}
	
If $m_H^F(\Sigma)<0$, then using the change of variables
\begin{align*}
- \sqrt{\frac{4 |\PW|_F}{|\Sigma|_F  }} m_H^F(\Sigma)e^{-t/2}= r^{\frac{p-1}{3-p}},
\end{align*}
we get by direct computation
\begin{align*}
\mathrm{Cap}_{F, p}(K) &\leq  \left( \frac{3-p}{p-1} \right)^{p-1} |\PW|_F^{\frac{p-1}{2}} |\Sigma|_F^{\frac{3-p}{2}} \left( \frac{\int_{\Sigma} H_F^2\,d\mu_F}{4|\PW|_F}  - 1\right)^{3-p}  \theta^{1-p},
\end{align*}
where
\begin{align*}
\theta := \int_{0}^{\left( \frac{\int_{\Sigma} H_F^2\,d\mu_F}{4|\PW|_F} -1 \right)^{\frac{3-p}{p-1}}} \left(1 +  r^{\frac{p-1}{3-p}}\right)^{-\frac{1}{2}}\,dr.
\end{align*}

	Should the equality hold, then \eqref{eq-inequality-surface} must be an equality. The H\"{o}lder inequality becomes an equality, which implies that $H_F$ is constant. So $\Sigma_t$ and then $\Sigma$ are translated scaled Wulff shapes (\cite{HLMG09}), which is impossible. Thus we cannot have the equality in this case. The proof is now complete.

\end{proof}

\begin{proof}[Proof of Corollary \ref{capf}]
Let $p = 2$ in Theorem~\ref{thm2}. If $m_H^F(\Sigma) = 0$, the conclusion follows immediately.

If $m_H^F(\Sigma) < 0$, we obtain
\begin{equation*}%\label{cor4ineq1}
\mathrm{Cap}_{F,2}\left( K\right) < \sqrt{|\PW|_F |\Sigma|_F} \left( \frac{1}{4|\PW|_F}\int_{\Sigma} H_F^2\,d\mu_F  - 1\right) \cdot \theta^{-1},
\end{equation*}
where
\begin{equation*}
		\theta = \int_{0}^{ \frac{\int_{\Sigma} H_F^2\,d\mu_F}{4|\PW|_F} -1 } \left(1 +  r\right)^{-\frac{1}{2}}\,dr = 2\left(\sqrt{ \frac{1}{4|\PW|_F}\int_{\Sigma} H_F^2\,d\mu_F}  - 1\right).
\end{equation*}
So we get
\begin{equation*}
\mathrm{Cap}_{F,2}( K) <  \frac{1}{2} \sqrt{|\PW|_F |\Sigma|_F} \left(1 + \sqrt{\frac{1}{4 |\PW|_F} \int_{\Sigma} H_F^2\,d\mu_F}\right).
\end{equation*}
Combining both cases we finish the proof of Corollary~\ref{capf}.
\end{proof}

\section{Proof of Theorem \ref{thm3}}\label{sec5}

\subsection{The case that $M=\pt K$ is convex.}
In this case we use the anisotropic unit speed flow
	\begin{equation*}
		\begin{cases}
X:N\times [0,T)\to \R^n,\\
			\pt_t X=\nu_F,\\
			X(N, 0)=M.
		\end{cases}
	\end{equation*}
	So we can derive $D \psi=F^{-1}(\nu)\nu$ and then
	\begin{align*}
		T_{p}(t)&=\int_{\psi(x)=t}\frac{F^p(D\psi)}{|D \psi|} d\mu=\int_{\psi(x)=t}F(\nu) d\mu.
\end{align*}
Note that on $M_t=\{x|\psi(x)=t\}$ we have
\begin{equation*}
X(z,t)=X(z,0)+t\nu_F(z,0),\quad z\in N.
\end{equation*}
Let $\Omega$ be the volume form on $\R^n$ and choose a local coordinate system $\{y^i\}_{i=1}^{n-1}$ on $N$ such that $\{\pt X(z,0)/\pt y^i\}_{i=1}^{n-1}$ correspond to the anisotropic principal directions of $M=X(N,0)$ at the point $X(z,0)$. Then we get
\begin{align*}
d\mu(X(z,t))&=\Omega\left(\nu,\frac{\pt X(z,t)}{\pt y^1},\dots,\frac{\pt X(z,t)}{\pt y^{n-1}}\right)dy^1\wedge \cdots \wedge dy^{n-1}\\
&=\left(\prod_{i=1}^{n-1}(1+\kappa_i^F t)\right)\Omega\left(\nu,\frac{\pt X(z,0)}{\pt y^1},\dots,\frac{\pt X(z,0)}{\pt y^{n-1}}\right)dy^1\wedge \cdots \wedge dy^{n-1}\\
&=\left(\prod_{i=1}^{n-1}(1+\kappa_i^F t)\right)d\mu(X(z,0)).
\end{align*}
As a consequence, we obtain (note that at $X(z,0)$ and $X(z,t)$ the outward unit normal is the same $\nu$)
\begin{align*}
		T_{p}(t)&=\int_{M}\prod_{i=1}^{n-1}(1+\kappa_i^F t)F(\nu)d\mu=\sum_{i=0}^{n-1}\int_{M} \sigma_i(\kappa^F) d\mu_F \cdot t^i.
	\end{align*}
	Therefore we get
	\begin{align*}
		\mathrm{Cap}_{F, p}(K)& \leq \left(\int_0^\infty  T_{p}^{-\frac{1}{p-1}}(t)dt\right)^{1-p}\\
&= \left(\int_0^\infty \left(\sum_{i=0}^{n-1}\int_{M} \sigma_i(\kappa^F) d\mu_F \cdot t^i\right)^{\frac{1}{1-p}}dt\right)^{1-p} .
	\end{align*}

Next assume that the equality holds. Then the hypersurfaces $M_t$ ($t\geq 0$) constructed above are exactly the level sets of the $p$-capacitary potential $u$ of the set $K$ and $\lambda(\psi(x))=u(x)$. In particular, on $M$ we obtain
\begin{align*}
F(Du)=|\lambda'(0)|F(D\psi)=|\lambda'(0)|F(F^{-1}(\nu)\nu)=|\lambda'(0)|,
\end{align*}
a constant. Then by \cite[Theorem~1.2]{BC18}, we conclude that $M$ is a translated scaled Wulff shape. So we finish the proof in this case.

	\subsection{The case that $M$ is star-shaped with respect to the origin.}
	In this case we expand $M$ by homothety so that $M_t=(1+t)M$ for $t\geq 0$. So we have
	\begin{equation*}
		\psi((1+t)X)=t,\quad X\in M.
	\end{equation*}
	Let $X_t:=(1+t)X$. So taking the derivative with respect to $t$ yields
	\begin{equation*}
		\langle D \psi(X_t),X\rangle=1.
	\end{equation*}
	Note that at points $X_t$ and $X$ the outward unit normal is the same $\nu$, and $D \psi(X_t)=|D \psi(X_t)|\nu$. Hence we get
	\begin{equation*}
		|D \psi(X_t)|=\frac{1}{\langle X,\nu\rangle}.
	\end{equation*}
	Then we have (note that $d\mu(X_t)=(1+t)^{n-1}d\mu(X)$)
	\begin{align*}
		T_{p}(t)&=\int_{M_t}\frac{F^p(D\psi)}{|D\psi|}d\mu=\int_{M_t}\frac{F^p(\nu)}{\langle X,\nu\rangle^{p-1}}d\mu=(1+t)^{n-1}\int_{M} h_F^{1-p}d\mu_F,
	\end{align*}
	where $h_F=\langle X,\nu\rangle /F(\nu)$ is the anisotropic support function. So we get
	\begin{align*}
		\mathrm{Cap}_{F, p}(K) & \leq \left(\int_0^\infty  T_{p}^{-\frac{1}{p-1}}(t)dt\right)^{1-p}\\
&= \left(\frac{n-p}{p-1}\right)^{p-1}\int_{M} h_F^{1-p}d\mu_F.
	\end{align*}

Next assume that the equality holds. Then again the hypersurfaces $M_t$ ($t\geq 0$) constructed above are the level sets of the $p$-capacitary potential $u$ of the set $K$. Let $\rho(t)=u(X)$ for $X\in M_t$, $t\geq 0$.

Now consider any two points $X,Y\in M$. Then $(1+t)X$ and $(1+t)Y$ belong to $M_t$ for $t\geq 0$. We apply Proposition~3.1 in \cite{XY20} (cf. Lemma~5.2 in \cite{AGHLV17}) to deduce
\begin{align*}
\mathrm{Cap}_{F,p}(K)^{1/(p-1)}&=\lim_{t\to +\infty}\frac{u((1+t)X)}{c({n,p})(F^0((1+t)X))^{(p-n)/(p-1)}}\\
&=(F^0(X))^{(n-p)/(p-1)}\lim_{t\to +\infty}\frac{\rho(t)}{c({n,p})(1+t)^{(p-n)/(p-1)}},
\end{align*}
where $c(n,p)$ is a constant depending only on $n$ and $p$. The same holds for $Y$. Thus we have $F^0(X)=F^0(Y)$, which means that $M$ is a scaled Wulff shape centered at the origin. So we finish the proof in this case as well.

\section{Proof of Theorem \ref{thm4}}\label{sec6}

\begin{proof}
%For a compact convex set $K$ with non-empty interior, we construct a family of compact convex sets
%\begin{align*}
%K_\varepsilon:=K+\varepsilon \overline{B},\quad \varepsilon>0,
%\end{align*}
%where $B$ is the open unit ball in $\R^n$. Note that $K_\varepsilon$ has a $C^{1,1}$ boundary. If we can prove Theorem~\ref{thm4} for $K_\varepsilon$, then taking $\varepsilon \to 0+$ and using Proposition~\ref{prop-property}, we get the conclusion for $K$ itself. So in the proof below we assume that $K$ has non-empty interior and a $C^{1,1}$ boundary.
%
%Let $u$ be the anisotropic $p$-capacitary potential for $K$, i.e., $u$ solves
%\begin{equation*}
%		\begin{cases}
%			\Delta_{F, p} u = 0 &\text{ in }\R^n\setminus K, \\
%			u(x) = 1 &\text{ on } \partial K,\\
%			u(x)\rightarrow  0&\text{ as } |x| \to \infty.
%		\end{cases}
%	\end{equation*}
%Note that $K$ is a compact convex set with non-empty interior and with a $C^{1,1}$ boundary. By \cite[Appendix~B]{BC18} or \cite[Lemma~4.3]{AGHLV17}, we know $Du\neq 0$ on $\R^n\setminus K$. Thus $u$ is smooth on $\R^n\setminus K$ with smooth level sets $\{x|u(x)=t\}$ for $0<t< 1$. Moreover, by \cite[Lemma~4.4]{AGHLV17}, the level set $\{x|u(x)= t\}$ is convex for $0<t< 1$. We will use these facts below.

Let $u$ be the anisotropic $p$-capacitary potential for $K$, i.e., $u$ solves \eqref{eq-capacitary-potential}. Moreover, thanks to Proposition~\ref{prop-capacitary-potential}, we know that the set $K_t$ := $\{x | u(x) \ge t\}$ is convex with smooth boundary for any $0<t<1$.

Now recall that for a convex body $K$ (i.e., a compact convex set with non-empty interior), its anisotropic perimeter $|\pt K|_F$ can be expressed as a mixed volume of $K$ and $\overline{\mathcal{W}}$. More precisely, we have
\begin{equation}\label{eq-mixed-volume}
|\pt K|_F=\frac{1}{n}V_{(1)}(K,\overline{\mathcal{W}}),
\end{equation}
where $V_{(i)}(K,\overline{\mathcal{W}})$ ($0\leq i\leq n$) is the $i$th mixed volume of the convex bodies $K$ and $\overline{\mathcal{W}}$ defined as in the expression
\begin{align*}
|K+t\overline{\mathcal{W}}|=\sum_{i=0}^n\binom{n}{i}t^iV_{(i)}(K,\overline{\mathcal{W}}),\quad t\geq 0.
\end{align*}
See the classical book \cite{Sch14} for the definition and properties on mixed volumes. In particular, it follows from \cite[(5.25)]{Sch14} that for any two convex bodies $K_1$ and $K_2$ with $K_1\subset K_2$ we have $V_{(i)}(K_1,\overline{\mathcal{W}})\leq V_{(i)}(K_2,\overline{\mathcal{W}})$, $0\leq i\leq n$.

In our case we know that $K$ and $K_t$ ($0<t<1$) are convex bodies with $K\subset K_t$, which implies $|\partial K|_F \leq |\pt K_t|_F$ for any $t\in (0,1)$. Furthermore, we make a key observation that in fact,
\begin{align}\label{eq-key-observation}
|\partial K|_F < |\pt K_t|_F,\quad \forall t\in (0,1).
\end{align}
This observation will be used in the argument for the equality case at the end of the proof and can be checked based on the fact that the Hausdorff distance of $K$ and $K_t$ is positive for a fixed $t\in (0,1)$.

Next applying the H\"{o}lder inequality, the co-area formula and the fact
\begin{align*}
\mathrm{Cap}_{F, p}(K) = \int_{ u(x) = t} \frac{F^p(Du)}{|Du|} \,d\mu \text{ for any } t \in (0, 1),
\end{align*}
which can be proved straightforwardly (see Lemma~4.2 in \cite{AGHLV17}; cf. Lemma~2.16 in \cite{AKP15}), we derive (note that $\nu=-Du/|Du| $ on $\{x | u(x) = t\}$)
\begin{align*}
|\partial K|_F &<  \int_{     u(x) = t} \left(|Du|^{\frac{p-1}{p}} F(\nu)\right) |Du|^{\frac{1-p}{p}}\,d\mu\\
			&\le \left( \int_{  u(x) = t }\frac{F^p(Du)}{|Du|}\,d\mu   \right)^{\frac{1}{p}}  \left( \int_{  u(x) = t } |Du|^{-1} \,d\mu \right)^{\frac{p-1}{p}}\\
			&= \mathrm{Cap}_{F, p}(K)^{\frac{1}{p}} \left( -\frac{d}{dt} | K_t| \right)^{\frac{p-1}{p}}.
\end{align*}
It follows that
\begin{equation*}
	\left(  \frac{|\partial K|_F}{\mathrm{Cap}_{F, p}(K)^{\frac{1}{p}}}\right)^{\frac{p}{p-1}} < - \frac{d}{dt} | K_t|.
\end{equation*}
Integrating the above inequality over $(t, 1)$ yields
\begin{equation}\label{eq-6.3}
	(1-t) \left(  \frac{|\partial K|_F}{\mathrm{Cap}_{F, p}(K)^{\frac{1}{p}}}\right)^{\frac{p}{p-1}} < |K_t| - |K|.
\end{equation}

%Now we applies to \cite{Mazya11}  2.2 Cor 2, which states
%
%\begin{equation}\label{cap_pf lower bound}
%	Cap_{F, p}(K) \ge \left(\int_{|K|}^{\infty} \frac{d\rho}{\left( \mathscr{C}(\rho) \right)^{\frac{p}{p-1}}}\right)^{1-p}
%\end{equation}
%
%Where $\mathscr{C}(\rho) := \inf\{|\partial K| |\;\text{for all}\; K\; \text{such that}\; |K| \ge \rho\}$.
%
%Use \cite{IS91} Cor 2.8 (or arxiv:1604.04302v4), we know anisotropic isoperimetric inequality
%
%\begin{equation*}
%	|\partial K|_F \ge n |\W|^{\frac{1}{n}} |K|^{\frac{n-1}{n}}
%\end{equation*}
%
%Gives $\mathscr{C}(\rho) = n |\W|^{\frac{1}{n}} \rho^{\frac{n-1}{n}}$, put in (\ref{cap_pf lower bound}), we obtained
Now we recall the following isocapacitary inequality (see, e.g., (2.2.8) in \cite{Maz11})
\begin{equation*}
	\begin{aligned}
		\mathrm{Cap}_{F, p}(K) \ge n |\W|^{\frac{p}{n}} \left(\frac{n-p}{p-1}\right)^{p-1} |K|^{\frac{n-p}{n}}.
	\end{aligned}
\end{equation*}
Applying the above inequality to $K_t$ gives us
\begin{equation*}
	|K_t| \le n^{\frac{-n}{n-p}} |\W|^{\frac{-p}{n-p}} \left( \frac{\mathrm{Cap}_{F, p}(K_t)}{\left((n-p)/(p-1)\right)^{p-1}} \right)^{\frac{n}{n-p}}.
\end{equation*}
Combining it with \eqref{eq-6.3}, we get
\begin{equation*}
	(1-t) \left(  \frac{|\partial K|_F}{\mathrm{Cap}_{F, p}(K)^{\frac{1}{p}}}\right)^{\frac{p}{p-1}} < n^{\frac{-n}{n-p}} |\W|^{\frac{-p}{n-p}} \left( \frac{\mathrm{Cap}_{F, p}(K_t)}{\left((n-p)/(p-1)\right)^{p-1}} \right)^{\frac{n}{n-p}} - |K|.
\end{equation*}
Notice that $u/t$ is the anisotropic $p$-capacitary potential for the set $K_t$. So we have
\begin{align*}
		\mathrm{Cap}_{F, p}\left(K_t\right) &= \int_{  u(x) = t } \frac{F^p(D(u/t))}{|D(u/t)|}\,d\mu= t^{1-p} \mathrm{Cap}_{F, p}(K),\quad 0<t<1.
\end{align*}
Then we get
\begin{align}
	&(1-t) \left(  \frac{|\partial K|_F}{\mathrm{Cap}_{F, p}(K)^{\frac{1}{p}}}\right)^{\frac{p}{p-1}}\nonumber \\
  <n^{\frac{-n}{n-p}} &|\W|^{\frac{-p}{n-p}} \left( \frac{t^{1-p} \mathrm{Cap}_{F, p}(K)}{\left((n-p)/(p-1)\right)^{p-1}} \right)^{\frac{n}{n-p}} - |K|,\quad t\in (0,1).\label{eq-key-inequality}
\end{align}

Next define
\begin{align*}
\varphi(t) :=& (1-t) \left(  \frac{|\partial K|_F}{\mathrm{Cap}_{F, p}(K)^{\frac{1}{p}}}\right)^{\frac{p}{p-1}} \\
& {}- n^{\frac{-n}{n-p}} |\W|^{\frac{-p}{n-p}} \left( \frac{t^{1-p} \mathrm{Cap}_{F, p}(K)}{\left((n-p)/(p-1)\right)^{p-1}} \right)^{\frac{n}{n-p}} + |K|,\\
\bar{t} :=& \frac{p-1}{n-p}n^{\frac{1}{1-p}} |\W|^{\frac{-p}{n(p-1)}}  |K|^{-\frac{n-p}{n(p-1)}} \mathrm{Cap}_{F, p}(K)^{\frac{1}{p-1}}.
\end{align*}
Note that by the isocapacitary inequality, we have $\bar{t}\geq 1$.

We claim that $\varphi(t) \leq  0$ on the interval $(0,\bar{t}]$ and $\varphi(t)<0$ on $(0,1)$. If $t \in (0, 1)$, the claim follows from  \eqref{eq-key-inequality}. If $1\leq t\leq \bar{t}$, it follows from the fact
\begin{align*}
&	(1-t) \left(  \frac{|\partial K|_F}{\mathrm{Cap}_{F, p}(K)^{\frac{1}{p}}}\right)^{\frac{p}{p-1}} \le 0\\
\le n^{\frac{-n}{n-p}} &|\W|^{\frac{-p}{n-p}} \left( \frac{t^{1-p} \mathrm{Cap}_{F, p}(K)}{\left((n-p)/(p-1)\right)^{p-1}} \right)^{\frac{n}{n-p}} - |K|.
\end{align*}
So we obtain the claim. By the claim we have $\sup_{t\in (0,\bar{t}]}\varphi(t)\leq 0$.

Next we compute
\begin{align*}
	\varphi'(t) =&  - \left(  \frac{|\partial K|_F}{\mathrm{Cap}_{F, p}(K)^{\frac{1}{p}}}\right)^{\frac{p}{p-1}}\\
 &{}- n^{\frac{-n}{n-p}} |\W|^{\frac{-p}{n-p}} \frac{n(1-p)}{n-p} \left( \frac{ \mathrm{Cap}_{F, p}(K)}{\left((n-p)/(p-1)\right)^{p-1}} \right)^{\frac{n}{n-p}} t^{-\frac{p(n-1)}{n-p}}.
\end{align*}
So the unique critical point $t_0$ of $\varphi(t)$ reads
\begin{equation*}
	t_0 = \left(n |\W|\right)^{-\frac{1}{n-1}} \frac{p-1}{n-p}  |\partial K|_F^{-\frac{n-p}{(n-1)(p-1)}} \mathrm{Cap}_{F, p}(K)^{\frac{1}{p-1}}.
\end{equation*}
Now recall the anisotropic isoperimetric inequality (see e.g. Theorem~20.8 in \cite{Mag12} or Corollary~2.8 in \cite{FM91})
\begin{align}\label{eq-isoperimetric}
|\pt K|_F\geq n|\mathcal{W}|^{1/n}|K|^{(n-1)/n},
\end{align}
with the equality holding if and only if $\pt K$ is a translated scaled Wulff shape. Using it we know $t_0 \leq \bar{t}$. Note that $\varphi(0^+)=-\infty$, $\varphi(t)$ is non-decreasing on $(0,t_0]$ and non-increasing on $[t_0,\bar{t}]$. So $\varphi(t_0)\leq 0$, which by direct computation is equivalent to
\begin{align*}
\frac{p(n-1)}{n(n-p)}|\partial K|_F^{\frac{n}{n-1}}& |\pt \W|_F^{\frac{1}{1-n}} \geq  |K|+ \frac{|\partial K|_F^{p/(p-1)}}{\mathrm{Cap}_{F, p}(K)^{1/(p-1)}},
\end{align*}
the desired inequality.

Last, if $\pt K$ is a translated scaled Wulff shape, we can check directly that the equality in the above inequality holds. Conversely, assume that we have the equality $\varphi(t_0)=0$, i.e.,
\begin{align*}
\frac{p(n-1)}{n(n-p)}|\partial K|_F^{\frac{n}{n-1}}& |\pt \W|_F^{\frac{1}{1-n}}= |K|+ \frac{|\partial K|_F^{p/(p-1)}}{\mathrm{Cap}_{F, p}(K)^{1/(p-1)}}.
\end{align*}
Using the anisotropic isoperimetric inequality \eqref{eq-isoperimetric} to replace $|K|$ leads to
\begin{align*}
\frac{p-1}{n-p}|\pt \mathcal{W}|^{-\frac{1}{n-1}}|\pt K|_F^{-\frac{n-p}{(p-1)(n-1)}}\mathrm{Cap}_{F,p}(K)^{\frac{1}{p-1}}\leq 1.
\end{align*}
This is nothing but $t_0\leq 1$. However, $t_0<1$ can not occur, since $\varphi(t_0)=0$ and by the previous claim we know $\varphi(t)<0$ for $t\in (0,1)$. (This is where we use our key observation \eqref{eq-key-observation}.) Thus we must have $t_0=1$ and then we conclude that $\pt K$ is a translated scaled Wulff shape from the rigidity part of the anisotropic isoperimetric inequality.

Now the proof of Theorem~\ref{thm4} is complete.
\end{proof}

\bibliographystyle{Plain}

\begin{thebibliography}{99}

%\bibitem{And14} Ben~Andrews, \emph{Moduli of continuity, isoperimetric profiles, and multi-point estimates in geometric heat equations}, Surveys in Differential Geometry \textbf{19} (2014), 1--47.


\bibitem{AFM20} Virginia~Agostiniani, Mattia~Fogagnolo, and Lorenzo~Mazzieri, \emph{Sharp geometric inequalities for closed hypersurfaces in manifolds with nonnegative Ricci curvature}, Invent. Math. \textbf{222} (2020), no.~3, 1033--1101.

\bibitem{AM15} Virginia~Agostiniani and Lorenzo~Mazzieri, \emph{Riemannian aspects of potential theory}, J. Math. Pures Appl. (9) \textbf{104} (2015), no.~3, 561--586.

\bibitem{AM20} Virginia~Agostiniani and Lorenzo~Mazzieri, \emph{Monotonicity formulas in potential theory}, Calc. Var. Partial Differential Equations \textbf{59} (2020), no.~1, Paper No.~6, 32 pp.

\bibitem{AGHLV17} M.~Akman, J.~Gong, J.~Hineman, J.~Lewis, and A.~Vogel, \emph{The Brunn--Minkowski inequality and a Minkowski problem for nonlinear capacity}, to appear in Mem. Amer. Math. Soc., arXiv:1709.00447.

\bibitem{And01} Ben~Andrews, \emph{Volume-preserving anisotropic mean curvature flow}, Indiana Univ. Math. J. \textbf{50} (2001), no.~2, 783--827.

\bibitem{BC18} Chiara~Bianchini and Giulio~Ciraolo, \emph{Wulff shape characterizations in overdetermined anisotropic elliptic problems}, Comm. Partial Differential Equations \textbf{43} (2018), no.~5, 790--820.

\bibitem{BCS16} Chiara~Bianchini, Giulio~Ciraolo, and Paolo~Salani, \emph{An overdetermined problem for the anisotropic capacity}, Calc. Var. Partial Differential Equations \textbf{55} (2016), no.~4, Art. 84, 24 pp.

\bibitem{BCS18} Chiara~Bianchini, Giulio~Ciraolo, and Paolo~Salani, \emph{Some overdetermined problems related to the anisotropic capacity}, J. Math. Anal. Appl. \textbf{465} (2018), no.~1, 211--219.

\bibitem{BM08} H.~Bray and P.~Miao, \emph{On the capacity of surfaces in manifolds with nonnegative scalar curvature}, Invent. Math. \textbf{172} (2008), no.~3, 459--475.

\bibitem{AKP15} A.~Colesanti, K.~Nystr\"{o}m, P.~Salani, J.~Xiao, D.~Yang, and G.~Zhang, \emph{The Hadamard variational formula and the Minkowski problem for $p$-capacity}, Adv. Math. \textbf{285} (2015), 1511--1588.

%\bibitem{Cru19} Tiarlos~Cruz, \emph{Capacity inequalities and rigidity of cornered/conical manifolds}, Ann. Global Anal. Geom. \textbf{55} (2019), no.~2, 281--298.

\bibitem{FM91} Irene~Fonseca and Stefan~M\"{u}ller, \emph{A uniqueness proof for the Wulff theorem}, Proc. Roy. Soc. Edinburgh Sect. A \textbf{119} (1991), no.~1-2, 125--136.

\bibitem{FS14} Alexandre~Freire and Fernando~Schwartz, \emph{Mass-capacity inequalities for conformally flat manifolds with boundary}, Comm. Partial Differential Equations \textbf{39} (2014), no.~1, 98--119.

\bibitem{Ger90} Claus~Gerhardt, \emph{Flow of nonconvex hypersurfaces into spheres}, J. Differential Geom. \textbf{32} (1990), no.~1, 299--314.

\bibitem{Gero73} Robert~Geroch, \emph{Energy extraction}, Annals of the New York Academy of Sciences \textbf{224} (1973), 108--117.

\bibitem{Haw68} Stephen~W.~Hawking, \emph{Gravitational radiation in an expanding universe}, Journal of Mathematical Physics \textbf{9} (1968), 598--604.

\bibitem{HLMG09} Yijun~He, Haizhong~Li, Hui~Ma, and Jianquan~Ge, \emph{Compact embedded hypersurfaces with constant higher order anisotropic mean curvatures}, Indiana Univ. Math. J. \textbf{58} (2009), no.~2, 853--868.

\bibitem{HKM06} Juha~Heinonen, Tero~Kilpel\"{a}inen, and Olli~Martio, \emph{Nonlinear potential theory of degenerate elliptic equations}, Unabridged republication of the 1993 original, Dover Publications, Inc., Mineola, NY, 2006.



\bibitem{HI01} Gerhard~Huisken and Tom~Ilmanen, \emph{The inverse mean curvature flow and the Riemannian Penrose inequality}, J. Differential Geom. \textbf{59} (2001), no.~3, 353--437.

\bibitem{LL01} Elliott~H.~Lieb and Michael~Loss, \emph{Analysis}, Second edition, Graduate Studies in Mathematics, \textbf{14}, American Mathematical Society, Providence, RI, 2001.

\bibitem{Mag12} Francesco~Maggi, \emph{Sets of finite perimeter and geometric variational problems: An introduction to geometric measure theory}, Cambridge Studies in Advanced Mathematics, \textbf{135}, Cambridge University Press, Cambridge, 2012.

\bibitem{Maz11} Vladimir~Maz'ya, \emph{Sobolev spaces with applications to elliptic partial differential equations}, Second, revised and augmented edition, Grundlehren der Mathematischen Wissenschaften [Fundamental Principles of Mathematical Sciences], \textbf{342}, Springer, Heidelberg, 2011.

\bibitem{Pol47} G.~P\'{o}lya, \emph{Estimating electrostatic capacity}, Amer. Math. Monthly \textbf{54} (1947), 201--206.

\bibitem{PS51} G.~P\'{o}lya and G.~Szeg\"{o}, \emph{Isoperimetric Inequalities in Mathematical Physics}, Annals of Mathematics Studies, no.~27, Princeton University Press, Princeton, N. J., 1951.

\bibitem{Rei76} Robert C.~Reilly, \emph{The relative differential geometry of nonparametric hypersurfaces}, Duke Math. J. \textbf{43} (1976), no.~4, 705--721.

\bibitem{Sch14} Rolf Schneider,  \emph{Convex bodies: the Brunn--Minkowski theory}, Second expanded edition, Encyclopedia of Mathematics and its Applications, vol. 151, Cambridge University Press, Cambridge, 2014.

\bibitem{Shu} M.~Shubin, \emph{Capacity and its applications}, \path{https://citeseerx.ist.psu.edu/viewdoc/download?doi=10.1.1.140.1404&rep=rep1&type=pdf}.

\bibitem{Sze31} G.~Szeg\"{o}, \emph{\"{U}ber einige neue Extremaleigenschaften der Kugel}, Math. Z. \textbf{33} (1931), no.~1, 419--425.

%\bibitem{Sze47} G.~Szeg\"{o},

\bibitem{Urb90} John~I.~E.~Urbas, \emph{On the expansion of starshaped hypersurfaces by symmetric functions of their principal curvatures}, Math. Z. \textbf{205} (1990), no.~3, 355--372.

\bibitem{Xia13} Chao~Xia, \emph{On an anisotropic Minkowski problem}, Indiana Univ. Math. J. \textbf{62} (2013), no.~5, 1399--1430.

\bibitem{Xia17} Chao~Xia, \emph{Inverse anisotropic mean curvature flow and a Minkowski type inequality}, Adv. Math. \textbf{315} (2017), 102--129.

\bibitem{XY20} Chao~Xia and Jiabin~Yin, \emph{Anisotropic $p$-capacity and anisotropic Minkowski inequality}, arXiv:2012.13933.

\bibitem{Xiao16} J.~Xiao, \emph{The $p$-harmonic capacity of an asymptotically flat $3$-manifold with non-negative scalar curvature}, Ann. Henri Poincar\'{e} \textbf{17} (2016), no.~8, 2265--2283.

\bibitem{Xiao17} Jie~Xiao, \emph{$P$-capacity vs surface-area}, Adv. Math. \textbf{308} (2017), 1318--1336.


\end{thebibliography}

\end{document}